\documentclass[12pt, draftclsnofoot, onecolumn]{IEEEtran}

\usepackage{amsmath,amsfonts}
\usepackage{algorithmic}
\usepackage{array}
\usepackage[caption=false,font=normalsize,labelfont=sf,textfont=sf]{subfig}
\usepackage{textcomp}
\usepackage{stfloats}
\usepackage{url}
\usepackage{verbatim}
\usepackage{graphicx}
\def\BibTeX{{\rm B\kern-.05em{\sc i\kern-.025em b}\kern-.08em
    T\kern-.1667em\lower.7ex\hbox{E}\kern-.125emX}}
\usepackage{balance}

\usepackage[utf8]{inputenc} % allow utf-8 input
\usepackage[T1]{fontenc}    % use 8-bit T1 fonts
\usepackage{hyperref}       % hyperlinks
\usepackage{booktabs}       % professional-quality tables
\usepackage{amsfonts}       % blackboard math symbols
\usepackage{nicefrac}       % compact symbols for 1/2, etc.
\usepackage{microtype}      % microtypography
\usepackage{xcolor}         % colors
\usepackage{import}
\usepackage{mathtools, amssymb, amsthm}
\usepackage[ruled,vlined]{algorithm2e}
\usepackage{tikz}
\usepackage{xcolor}
\usepackage{siunitx}
\usepackage{enumerate}
\usepackage{comment}
\usepackage{enumitem}
\usepackage[backend=biber]{biblatex}
\usepackage[capitalize,noabbrev]{cleveref}
\usepackage{tabularx} 

\DeclareLanguageMapping{english}{english-apa}
\newcounter{commentcounter}

\counterwithin*{commentcounter}{section}

\theoremstyle{plain}
\newtheorem{theorem}{Theorem}[section]

\newtheorem{lemma}[theorem]{Lemma}

\theoremstyle{definition}
\newtheorem{definition}[theorem]{Definition}
\newtheorem{assumption}[theorem]{Assumption}
\theoremstyle{remark}
\newtheorem{remark}[theorem]{Remark}

\DeclareMathOperator*{\argmin}{arg\,min}

\newcommand{\prox}{\mathrm{Prox}}

\newcommand{\RR}{\mathbb{R}}
\newcommand{\EE}{\mathbb{E}}

\newcommand{\PP}{\mathbb{P}}
\newcommand{\CC}{\mathbb{C}}

\newcommand{\cC}{{\mathcal{C}}}

\newcommand{\cE}{{\mathcal{E}}}
\newcommand{\cF}{{\mathcal{F}}}
\newcommand{\cG}{{\mathcal{G}}}

\newcommand{\cL}{{\mathcal{L}}}

\newcommand{\cN}{{\mathcal{N}}}
\newcommand{\cO}{{\mathcal{O}}}

\newcommand{\cS}{{\mathcal{S}}}
\newcommand{\cT}{{\mathcal{T}}}

\newcommand{\cV}{{\mathcal{V}}}

\newcommand{\va}{{\boldsymbol{a}}}
\newcommand{\vb}{{\boldsymbol{b}}}
\newcommand{\vc}{{\boldsymbol{c}}}

\newcommand{\ve}{{\boldsymbol{e}}}

\newcommand{\vs}{{\boldsymbol{s}}}
\newcommand{\vt}{{\boldsymbol{t}}}
\newcommand{\vu}{{\boldsymbol{u}}}
\newcommand{\vv}{{\boldsymbol{v}}}

\newcommand{\vx}{{\boldsymbol{x}}}
\newcommand{\vy}{{\boldsymbol{y}}}
\newcommand{\vz}{{\boldsymbol{z}}}

\newcommand{\vI}{{\boldsymbol{I}}}

\newcommand{\vlambda}{{\boldsymbol{\lambda}}}

\newcommand{\vxi}{{\boldsymbol{\xi}}}

\author{Behnam Mafakheri$^{1}$, \IEEEmembership{IEEE
} Jonathan H. Manton$^{1}$, \IEEEmembership{IEEE} and Iman Shames$^{2}$, \IEEEmembership{IEEE}% 
\thanks{*This work was supported by the Australian Research Council under the Discovery Projects funding scheme (DP210102454).}% <-this % stops a space
\thanks{$^{1}$B. Mafakheri and J. H. Manton are with Department of Electrical and Electronic Engineering,
        University of Melbourne, VIC 3010, Australia
        {\tt\small \{mafakherib, jmanton\}@unimelb.edu.au}}%
\thanks{$^{2}$I. Shames is  with the CIICADA Lab, School of Engineering, The Australian National University, Canberra, ACT 2601, Australia
        {\tt\small iman.shames@anu.edu.au}}%
}

\begin{document}
\title{An Asynchronous Decentralised Optimisation Algorithm for Nonconvex Problems}

% \markboth{Journal of \LaTeX\ Class Files,~Vol.~18, No.~9, September~2020}%
% {How to Use the IEEEtran \LaTeX \ Templates}

\maketitle

\begin{abstract}
    In this paper, we consider nonconvex decentralised optimisation and learning over a network of distributed agents. We develop an ADMM algorithm based on the Randomised Block Coordinate Douglas-Rachford splitting method which enables agents in the network to distributedly and asynchronously compute a set of first-order stationary solutions of the problem. To the best of our knowledge, this is the first decentralised and asynchronous algorithm for solving nonconvex optimisation problems with convergence proof. The numerical examples demonstrate the efficiency of the proposed algorithm for distributed Phase Retrieval and sparse Principal Component Analysis problems. 
\end{abstract}


\section{Introduction}
\IEEEPARstart{I}n the distributed optimisation problem, a set of $n$ computing nodes/agents aim to collaboratively solve optimisation problems of the following archetypal form: 
\begin{align}\label{my_P}
    \begin{array}{ll}
        \underset{\vx \in \RR^d}{\mbox{minimize}} &\displaystyle{F(\vx):=\sum_{i=1}^n f_i(\vx) + g(\vx)},   
    \end{array}
\end{align}
It is assumed that each agent $i\in \{1, 2, \dots, n\}$ individually possesses access to the local objective function $f_i(\vx)$ and the function $g(\vx)$. {This problem is a crucial element of distributed decision-making, with diverse applications including but not limited to compressed sensing, dictionary learning, power allocation in wireless ad hoc networks, distributed clustering, and empirical risk minimisation.} We next provide details of some examples of large scale networks and applications that arise within such settings.

\emph{Example 1} LASSO (Least Absolute Shrinkage and Selection Operator) is the problem where $F(x) = \sum_{i=1}^n \|A_i\vx - \vb_i \|^2 + \lambda \|\vx\|_1$ for $A_i \in \RR^{m_i\times d}$, $\vb_i \in \RR^{m_i}$, and $\lambda >0$. It is the ordinary least square problem with an $\ell_1$ penalisation in which each agent $i$ has access to $m_i$ measurements. 
LASSO is also one of the main models of compressed sensing \cite{donoho2006compressed}. 

\emph {Example 2} Phase Retrieval is the problem of detecting a signal from noisy measurements which capture the square of magnitude of the Fourier transform of the signal \cite{jaganathan2016phase, dong2023phase}. In other words, we wish to recover the phase of the signal with magnitude measurements. In a distributed setting, the problem can be written in the general form of $F(\vx) = \sum_{i=1}^n f_i(\vx)$ where $f_i(\vx) = \frac{1}{m_i} \sum_{j=1}^{m_i} (b_{ij} - |\langle \vx, \vt_{ij} \rangle |^2)^2$, and each agent $i$ has $m_i$ measurements $(b_{ij}, \vt_{ij})_{j=1}^{m_i}$.

{\emph{Example 3} Sparse Principal Component Analysis (SPCA) is a dimensionality reduction and machine learning method used to simplify a large data set into a smaller set while still maintaining significant patterns and trends and finds extensive applications in science and engineering, as demonstrated, for example, in \cite{guan2009sparse, sigg2008expectation}. We assume that agent $i$ has only access to a part of the data set shown by the matrix $P_i$, referred to as \emph{mini-batch} in the literature. One can formulate finding the first sparse principal component as the following distributed optimisation problem
$        \min_{\vx\in \RR^p} \sum_{i=1}^n -\|P_i\vx \|^2 + g(\vx), \text{subject to } \| \vx\|^2 \leq 1.$
Again, $g$ is a sparsity-promoting regulariser function which can be in the form of $\ell_0$ norm, $\ell_1$ norm, and the log sum penalty (LSP).}

\emph{Example 4} ERM (Empirical Risk Minimisation) is the specific context in which the vector $\vx$ is interpreted as the features to be learnt, with $f_i$ representing the cost function associated with a minibatch of data, and $g$ being used to impose desired structures on the solution (e.g., $\ell_1$ norm for sparsity) and/or enforce certain constraints. In machine learning applications, this serves to control the complexity of the model and prevent overfitting during the training process \cite[Ch. ~7]{Goodfellow-et-al-2016}. This exemplifies the majority of current deep learning optimisation problems.

In Example 1, it is evident that the local loss functions $f_i$ are convex, whereas in Examples 2-4, they are nonconvex. Another significant concern inherent in distributed computations is the emergence of asynchrony, wherein nodes are anticipated to execute their local computations at different speeds. In this work, we focus on the nonconvex problems in asynchronous setting. {More specifically, we consider ADMM-based algorithms for the distributed optimisation problem \eqref{my_P} over an arbitrary connected underlying communication graph. It is noteworthy that even in the synchronous setting, the literature for ADMM-based distributed optimisation algorithms on arbitrary connected graphs is sparse, see for example \cite{hajinezhad2017distributed} and \cite{yi2022sublinear}. We have discussed this in details in related works section.}

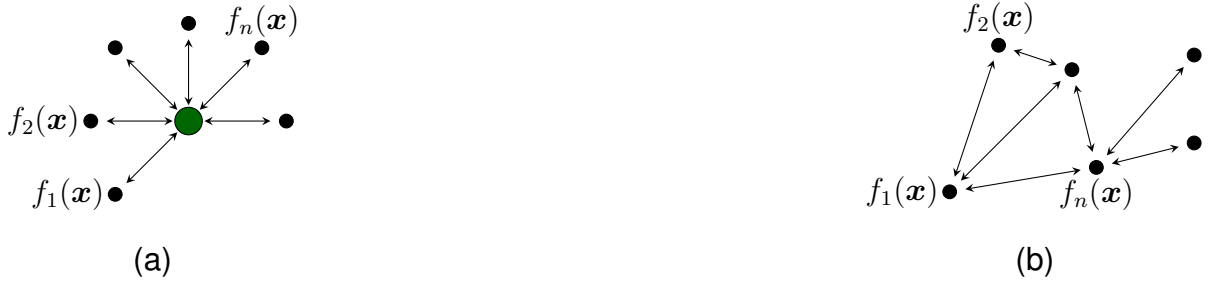
\begin{figure}[!t]
    \centering
    % \subfloat[b]{0.45\textwidth}
    \subfloat[]{
        \begin{tikzpicture}[scale = .65,
        trans/.style={<->,shorten >=6pt,shorten <=6pt,>=stealth}]
            \filldraw[black] (0.5,0.5) circle (4pt) node[anchor=east]{\color{black}{$f_1(\boldsymbol{x})$}};
            \filldraw[fill=black!60!green] (2, 2) circle (8pt) node[anchor=south]{\color{black}};
            \filldraw[black] (3.5, 3.5) circle (4pt) node[anchor=south]{\color{black}{$f_n(\boldsymbol{x})$}};
            \filldraw[black] (4, 2) circle (4pt) node[anchor=west]{\color{black}};
            \filldraw[black] (2, 4) circle (4pt) 
            node[anchor=west]{\color{black}};
            \filldraw[black] (0, 2) circle (4pt) node[anchor=east]{\color{black}{$f_2(\boldsymbol{x})$}};

            \filldraw[black] (0.5, 3.5) circle (4pt) node[anchor=west]{\color{black}};
            
            \draw[trans] (2,2) -- (0.5,0.5);
            \draw[trans] (2,2) -- (0, 2);
            \draw[trans] (2,2) -- (3.5,3.5);
            \draw[trans] (2,2) -- (4,2);
            \draw[trans] (2,2) -- (2, 4);
            \draw[trans] (2,2) -- (0.5,3.5);
        \end{tikzpicture}
        \label{fig:star_net}}
        % \caption{Distributed but not decentralised, nodes exchange information with the central node} 
    \hfill
    % \subfloat[b]{0.45\textwidth}
    \subfloat[]{  
            \begin{tikzpicture}[scale = .65,
        trans/.style={<->,shorten >=6pt,shorten <=6pt,>=stealth}]
            \filldraw[black] (0,0) circle (4pt) node[anchor=east]{\color{black}{$f_1(\boldsymbol{x})$}};
            \filldraw[black] (1, 3) circle (4pt) node[anchor=south]{\color{black}{$f_2(\boldsymbol{x})$}};
            \filldraw[black] (3, 0.5) circle (4pt) node[anchor=north]{\color{black}{$f_n(\boldsymbol{x})$}};
            \filldraw[black] (2.5, 2.5) circle (4pt) node[anchor=west]{\color{black}};
            \filldraw[black] (5, 1) circle (4pt) 
            node[anchor=west]{\color{black}};
            \filldraw[black] (5, 2.8) circle (4pt) node[anchor=west]{\color{black}};
            
            \draw[trans] (0,0) -- (1,3);
            \draw[trans] (2.5, 2.5) --  (3,0.5);
            \draw[trans] (1,3) -- (2.5, 2.5);
            \draw[trans] (3,0.5) -- (5,1);
            \draw[trans] (0,0) -- (3,0.5);
            \draw[trans] (0,0) -- (2.5,2.5);
            \draw[trans] (3,0.5) -- (5,2.8);
        \end{tikzpicture}
        \label{fig:mesh_net}}
            % \caption{Decentralised, nodes can exchange information with their immediate neighbours}
            
    % \end{subfigure}
    \caption{(a) Distributed and Centralised communication graph, (b) Distributed and Decentralised communication graph}
    \label{Dis-Dec}
\end{figure}

The traditional optimisation methodology requires access to all problem data in a central computing node. In the case of large-scale problems, centralised processing and storage present formidable challenges. Furthermore, in systems where agents are geographically distributed, aggregating all data in a central node is either inefficient or impractical due to energy constraints, link/hardware failures, and privacy concerns \cite{nedic2018multi}. Distributed/parallel optimisation typically involves the presence of a central node in the network that collects necessary information from nodes, updates the optimisation variable, and broadcasts it to all nodes (see Figure \ref{fig:star_net}). While this architecture shares some drawbacks with classical centralised algorithms, such as privacy concerns, it also introduces a single point of failure, which is undesirable in certain applications. Motivated by these practical challenges, our goal is to optimise in a decentralised architecture, which means that the architecture lacks a central node, and the nodes can only communicate (send or receive information) with their adjacent nodes in an underlying communication graph (see Figure \ref{fig:mesh_net}). In this paper, we make a distinction between \emph{distributed} and \emph{decentralised} methods, where decentralised methods specifically assume the absence of a central node.

% \begin{figure*}[!t]
% \centering
% \subfloat[]{\includegraphics[width=2.5in]{fig1}%
% \label{fig_first_case}}
% \hfil
% \subfloat[]{\includegraphics[width=2.5in]{fig1}%
% \label{fig_second_case}}
% \caption{Dae. Ad quatur autat ut porepel itemoles dolor autem fuga. Bus quia con nessunti as remo di quatus non perum que nimus. (a) Case I. (b) Case II.}
% \label{fig_sim}
% \end{figure*}

The challenge of agent asynchrony naturally arises when addressing the distributed optimisation problem \eqref{my_P} in a distributed manner. Given the heterogeneity of real-world nodes in terms of processing power, quality of communication links, and energy availability, devised distributed algorithms must tolerate variations in computation and communication delays among nodes. This adaptability enables some nodes to compute more quickly and communicate more frequently. Introducing an asynchronous algorithm can offer benefits, such as reducing agents' idle time, lowering power consumption, and alleviating congestion in the communication network. It should be noted that the convergence of the asynchronous variant of an algorithm is not guaranteed \cite[p.~103]{bertsekas1989parallel}. The majority of existing algorithms for solving \eqref{my_P} require the communication network to adopt a star graph configuration, involving a central node, or operate under a \emph{synchronous} paradigm. We extend this approach to align with a decentralised asynchronous framework, where a message-passing model is employed. In this setting, each worker is allowed to perform local updates with partially updated information received from its neighbours, and agents can converge to a critical point of the problem.

\textbf{Related Works.} The literature for convex distributed consensus optimisation methods is vast, for a recent survey on distributed convex optimisation we refer the reader to \cite{yang2019survey}. 

Distributed algorithms for solving \eqref{my_P} in nonconvex setting belong to one of the following two classes. The first is based on the gradient descent method and its variants. These algorithms are based on a (sub)gradient descent step followed by an averaging with neighbours and can be extended to time-varying and/or directed underlying graphs \cite{nedic2015dist, tatarenko2017non}. Some variants of Federated Learning (FL), an emerging branch of distributed learning that was first introduced in \cite{mcmahan2017communication}, are in this class of methods. FedProx \cite{li2020federated} tackles the heterogeneity of nodes in federated networks and an asynchronous block coordinate scheme has been proposed in \cite{wu2021federated}. The second class is based on operator splitting methods applied to the primal or dual problems. In operator splitting-based methods, slack variables along with equality constraints are introduced to decouple the objective function. At each step, some Lagrangian-related functions are optimised for a fixed dual variable, then the dual variables are updated accordingly. These methods are preferred when agents can solve their local optimisation problem efficiently. In this work, we explore the Alternating Direction Method of Multipliers (ADMM), a powerful operator splitting-based method widely recognised for its effectiveness in numerically solving optimisation problems \cite{boyd2011distributed, deng2016global, hong2016convergence}. However, existing distributed ADMM algorithms face limitations in their applicability, as they either assume a communication network structured as a star graph with a central node or operate in a synchronous setting. In scenarios where a central node is present in the graph, previous works \cite{hong2016convergence, chang2016asynchronous, hong2017distributed} have shown that ADMM converges in an asynchronous setting to a stationary point of the problem with a sublinear rate. Modified ADMM versions tailored for arbitrary connected communication graphs have been proposed by authors in \cite{pmlr-v70-hong17a, yi2022sublinear, mafakheri2023distributed}. Specifically, the Prox-PDA algorithm in \cite{pmlr-v70-hong17a} achieves sublinear convergence to a stationary point of the problem \eqref{my_P}. In a similar problem context, \cite{yi2022sublinear} proves the convergence of a modified ADMM version to the global minimum under the condition that the global cost function satisfies the Polyak-{\L}ojasiewicz condition and \cite{mafakheri2023distributed} has considered the decentralised problem in the presence of nonsmooth convex regulariser function $g$. {In \cite{guo2017asynchronous} an ADMM based asynchronous algorithm is developed for Optimal Power Flow problem in power systems where each agent only waits until it receives updated information from a subset of its neighbours in the graph. Authors have proved the convergence of the algorithm under the assumptions of compactness of the constrained set, restricted prox-regularity of cost functions, and boundedness of the sequence of dual variable.} In \cite{latafat2022block} global exponential convergence to the minimiser of the problem has been established under Kurdyka-{\L}ojasiewicz for the Forward-Backward splitting method. Closely related to our work are the works \cite{tran2021feddr} and \cite{wang2022fedadmm} where authors look at the asynchronous versions of DRS and ADMM, respectively. In both cases, the proposed algorithm needs a central node to perform the aggregation steps at each iteration. {In this paper, we propose an asynchronous ADMM-based algorithms for arbitrary graphs similar to the synchronous algorithms in \cite{pmlr-v70-hong17a} and \cite{yi2022sublinear}. It is noteworthy that the proposed algorithm in this paper is not the same as asynchronous version of \cite{pmlr-v70-hong17a} and \cite{yi2022sublinear} since their convergence under asynchronous condition is not guaranteed. This is further discussed in the Numerical Experiments section. Note that even under star topology, extending a synchronous algorithm to its asynchronous version and proving the convergence is non-trivial, see \cite{hong2016convergence} and \cite{chang2016asynchronous} for example.} The class of problems along with the key features of the algorithms proposed in these papers are summarised in Table \ref{table: literature}.

To address these limitations, we propose a decentralised ADMM-based algorithm that operates asynchronously. This means that certain agents within the network can awaken at random intervals to update their local parameters, all without the need for a global clock or a central node.

\begin{table*}[tbp]
    
    \centering
    \caption{Comparison with state-of-art distributed asynchronous algorithms.}
    \begin{tabularx}{\textwidth}{|>{\centering\arraybackslash}X|>{\centering\arraybackslash}X|>{\centering\arraybackslash}X|>{\centering\arraybackslash}X|>
    {\centering\arraybackslash}X|>{\centering\arraybackslash}X|}
        \hline
        \centering
        ADMM-based Algorithm & Class of Objective Functions & Parallel without a Central Node & Solving Subproblems &Asynchrony & Convergence to KKT Points\\
        \hline
        \hline

        AD-ADMM \cite{chang2016asynchronous} & Nonconvex with a nonsmooth $g$ & & Exact and non-exact & Partial & Deterministic $\cO(1/\sqrt{T})$\\
        \hline
        Gradient-based ADMM \cite{hong2017distributed} & Nonconvex with a nonsmooth $g$ & & Non-exact, using gradient only  & Partial & Deterministic $\cO(1/\sqrt{T})$\\
        \hline
        Prox-PDA \cite{pmlr-v70-hong17a} & Nonconvex & & Exact & \checkmark &  $\cO(1/\sqrt{T}) $\\
        \hline
        ASYMM \cite{farina2019distributed} & Nonconvex & & Exact & \checkmark  & Deterministic  \\
        \hline 
        Modified ADMM \cite{yi2022sublinear} & Nonconvex & \checkmark & Exact & & $\cO(1/\sqrt{T})$ \\
        \hline 
        L-ADMM \cite{yi2022sublinear} & Nonconvex (P-\L) & \checkmark & Non-exact, using gradient only & & Linear \\
        \hline
        Decentralised ADMM \cite{mafakheri2023distributed} &  Nonconvex with a nonsmooth $g$ & \checkmark & Exact & & Asymptotic \\
        \hline
        FedDR \cite{tran2021feddr} & Nonconvex & & Non-exact & \checkmark & Deterministic $\cO(1/\sqrt{T})$ \\
        \hline 
        This Work & Nonconvex & \checkmark & Exact & Partial & almost surely \\
        \hline
        
    \end{tabularx}
    \label{table: literature}
\end{table*}

\textbf{Our Contribution.} The focus of this paper is addressing the issue of asynchrony in decentralised optimisation. Our contributions and novelty can be summarised as follows.
\begin{itemize}
    \item We demonstrate the equivalence of 
 a version of the asynchronous ADMM algorithm and Block Coordinate Douglas-Rachford method (Theorem \ref{thm: asynaddm-rbcdrs equiv}). Such new insights can be of independent interest to optimisation community for understanding the asynchrony in distributed algorithms.
    \item We develop a new decentralised asynchronous ADMM-based algorithm for problem \eqref{my_P}, by combining the Douglas-Rachford Splitting technique, randomised block coordinate methods, and using the established equivalence mentioned above (see Algorithm \ref{alg:async-admm}). We assume a fixed and undirected underlying communication graph and decompose it into subgraphs in which each subgroup of agents can wake up independently and perform their local updates.
    \item We prove that if nodes are activated often enough and the penalty parameter of the algorithm is chosen large enough, then the proposed algorithm converges to the set of stationary points of the problem, starting from any initial point, with probability one (See Theorem \ref{thm: main result} for formal statement of the theorem).
\end{itemize}

\textbf{Organisation}: 
The paper is organised as follows. 
In the next section, we provide some preliminaries including notations,
definitions and mathematical backgrounds. In Section \ref{sec: alg and results}
we present the algorithm and the main results. In Section \ref{sec: Algorithm},
the proposed algorithm is given. In Section \ref{sec: alg derivation} we give
the details of derivation of the algorithm and in Section \ref{sec: proof BC-DRS}
the equivalence of Block Coordinate DRS and asynchronous ADMM is formally established. This result provides the theoretical grounding for the proposed asynchronous algorithm. Also, we establish the global convergence of the
algorithm in this section. Numerical comparisons and concluding remarks are
given in the last two sections.

\section {Notations and Preliminaries} \label{notations and preliminaries}

\noindent Throughout this paper, we represent the set of real numbers and extended real numbers as $\RR$ and $\overline{\RR}:=\RR \cup \{+\infty\}$. The complex numbers are denoted by $\CC$, and for a complex vector $\vz\in \CC^p$, the real and imaginary parts are indicated by $\Re{(z)}$ and $\Im{(z)}$, respectively. The set $[n]$ refers to the positive integers $\{1, 2, \dots, n \}$. For vectors $\va, \vb\in \RR^p$, the Euclidean inner product and corresponding norms are expressed as $\langle \va , \vb \rangle$ and $\| \va\|$, respectively. We use $col(\vx_1, \vx_2, \dots, \vx_n)$ to represent the stacking of vectors $\vx_1, \vx_2, \dots, \vx_n$ on top of each other, creating $[\vx_1^T, \vx_2^T, \dots, \vx_n^T]^T$. For a set $\cS \subset \RR^p$ and a matrix $A\in \RR^{m\times p}$, we define $A\cS:=\{A\vs \mid \vs \in \cS\}$. The indicator function of set $\cS$, denoted by $\iota_\cS(\cdot)$, is zero for members of set $\cS$ and infinity elsewhere. The Kronecker product of two matrices $A$ and $B$ with arbitrary sizes is denoted by $A\otimes B$. The vector of all ones with the appropriate length is denoted by $\boldsymbol{1}$. By $[\va]_i$, we denote the $i$-th element of vector $\va$, and $\ve_i$ is the vector that is one at positions of the $i$-th block and zero elsewhere. A graph $\mathcal{G}$ is defined by its set of vertices and edges. We represent it as $\mathcal{G}=([n], \mathcal{E})$, where $[n]$ is the set of nodes and $\mathcal{E} \subseteq [n]\times [n]$ is the set of edges. The presence of edge $(i,j)\in \mathcal{E}$ indicates that nodes $i$ and $j$ are adjacent, and node $i$ can receive information from node $j$. In this study, we assume that the graph $\cG$ is undirected, meaning that if $(i,j)\in \cE$, then $(j,i)\in \cE$, and there is a path between every pair of nodes, i.e., the graph is connected.

\section{Algorithm Development} \label{sec: alg and results}

\subsection{The asynchronous algorithm}\label{sec: Algorithm}
We consider the problem of the following form
\begin{align}\label{eq: optim_problem}
    \begin{array}{ll}
        \underset{\vx \in \RR^d}{\mbox{minimize}}&\displaystyle{f(\vx):=\sum_{i=1}^n f_i(\vx)}. 
    \end{array}
\end{align}
We assume that $n$ agents can communicate through a graph $\cG$  where each node can communicate with its immediate neighbours in the graph. Let a collection of $m$ not necessarily disjoint subsets of $[n]$ denoted by $\cV_1,
\dots, \cV_m$ such that $\bigcup_{j=1}^m \cV_j = [n]$, induce subgraphs $\cG_1, \dots, \cG_m$ where $\cG_j =
(\cV_j, \cE_j)$ with corresponding edge sets $\cE_j$.  We assume the following assumption over graph $\cG$ and the subgraphs.
\begin{assumption} \label{assp: subgraphs}
    Graph $\cG$ and each sub-graph $\cG_j, j\in [m]$ are connected.
\end{assumption}

Let $\vx_{\cV_i} := \textrm{col}(\{ \vx_j\}_{j\in \cV_i}) \in \RR^{d|\mathcal{V}_i|}$, for all $i\in [m]$, where $\textrm{col}(\va_1, \dots, \va_\ell) := [\va_1^T, \dots, \va_\ell^T]^T$. Define $\cN_i := \{j \mid i \in  \cV_j\}$. The proposed method is given in Algorithm \ref{alg:async-admm} for which the details of derivation will be given in Section \ref{sec: alg derivation}.

\begin{algorithm}[tb]
   \caption{Asynchronous ADMM}
   \label{alg:async-admm}
   At each iteration $k$, draw $\cS^k\subseteq [m]$ uniformly at random ($\cS^k$ is the set of subgraphs that update at iteration $k$)
   
    \For{$i \in \cS^k$}
    {
       \For{$j \in \cV_i$}
       {
             
             $\tilde{\vx}_j = \frac{1}{|\cN_j|}\sum_{r\in\cN_j} (\vz_{r_j} - \frac{1}{\beta}\vy_{r_j})$,  
             $\vx_j^+ \in \argmin_{\vu} f_j(\vu) + \frac{\beta |\cN_j|}{2}\|\vu -\tilde{\vx_j} \|^2$
        }
        $\bar{\vz}_i^+  = \frac{1}{|\cV_i|} \sum_{j\in \cV_i} (\vx_j^+ + \frac{\vy_{i_j}}{\beta})$, 
        $\vz_i^+ = (\underbrace{\bar{\vz}_i^+, \dots, \bar{\vz}_i^+}_{|\cV_i|\ \text{times}})$
        
        \For{$j \in \cV_i$}
        {
            $\vy_{i_j}^+ = \vy_{i_j} + \beta(\vx_j^+ - \bar{\vz}_i^+).$
        }
    }   
    \For {$i \notin \cS^k$}
    {
        \For {$j \in \cV_i$} 
        {
            $\vx_{j}^+ = \vx_{j}.$
        }
        {$\vz_i^+ = \vz_i, \vy_i^+ = \vy_i$}

    }
    
    $k \leftarrow k+1$.
\end{algorithm}
\begin{remark}
    Algorithm \ref{alg:async-admm} operates asynchronously, wherein, during each iteration, a randomly selected subset of agents $\cV_i$ is activated to update their local variables, while the remaining agents remain inactive. In other words, all nodes in a subgraph will perform their updates synchronously while the variables associated with other subgraphs will not change and hence we call the algorithm, \emph{partially} asynchronous.
\end{remark}

\begin{remark}
    The choice of sub-graphs $\cG_i$ is arbitrary. As long as the graph is connected (Assumption \ref{assp: subgraphs}) and the union of all sub-graphs is connected, then the algorithm is still well defined and the main result (in Theorem \ref{thm: asynaddm-rbcdrs equiv}) holds. For example, one can choose $\cG_i$s such that $\cup_{i=1}^m \cG_i$ constructs the smallest tree or a ring in the graph $\cG$.
\end{remark}

\begin{figure}[ht]
            \centering
            \begin{tikzpicture}[scale = 0.65,
        trans/.style={<->,shorten >=6pt,shorten <=6pt,>=stealth}]

            % \draw[help lines] (-1,-1) grid (6,5);

            \draw[black] (0,0) circle (8pt) node{\color{black}{$1$}};
            \draw[black] (1, 3) circle (8pt) node{\color{black}{$2$}};
            \draw[black] (3, 0.5) circle (8pt) node{\color{black}{$3$}};
            \draw[black] (2.5, 2.5) circle (8pt) node{\color{black}{$4$}};
            \draw[black] (5, 1) circle (8pt) 
            node{\color{black}{$5$}};
            \draw[black] (5, 2.8) circle (8pt) node{\color{black}{$6$}};
            
            \draw[trans] (0,0) -- (1,3);
            \draw[trans] (2.5, 2.5) --  (3,0.5);
            \draw[trans] (1,3) -- (2.5, 2.5);
            \draw[trans] (3,0.5) -- (5,1);
            \draw[trans] (0,0) -- (3,0.5);
            \draw[trans] (0,0) -- (2.5,2.5);
            \draw[trans] (2.5,2.5) -- (5,2.8);

            \draw[dashed] (1.26,1.2) circle (2.15);           
            \draw[dashed] (3.8,2.6) ellipse (1.7 and .8);
            \draw[dashed] (4,.72) ellipse (1.5 and .8);
            \node[above] at (0,3) {$\cG_1$};
            \node[above] at (5,3.2) {$\cG_2$};
            \node[right] at (5.4,1) {$\cG_3$};

        \end{tikzpicture}
            \caption{An example with $m=3$ where $\cV_1 = \{1,2,3, 4\}, \cV_2 = \{3, 5 \}, \cV_3 = \{4, 6 \}.$}
            \label{fig: example graph}
\end{figure}
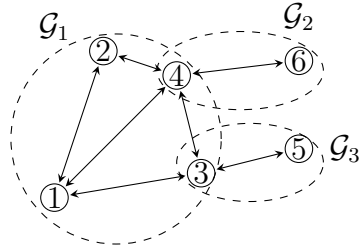

To give some insight into the formulation above, we give an example. Consider
the communication graph given in Fig \ref{fig: example graph}, where $n=6$,
$m=3$, $\vx = (\vx_1, \vx_2, \dots, \vx_6)$. We can write the problem as $\sum_{i=1}^6 f_i(\vx_i) + \iota_{\{\vx_1 = \vx_2 = \vx_3 = \vx_4\}}(\vx_1, \dots, \vx_4) + \iota_{\{\vx_3 = \vx_5\}}(\vx_3, \vx_5) + \iota_{\{\vx_4 = \vx_6\}}(\vx_4, \vx_6)$ where $\iota$ is the indicator function. Note that based on the
decomposition $\cV_1, \cV_2$ and $\cV_3$ given in the figure, the indicator
functions equals to zero if and only if
$\vx_1=\vx_2=\vx_3=\vx_4$, $\vx_3 = \vx_5$ and $\vx_4=\vx_6$. Thus, all $\vx_i$'s
are equal.

In the following section, we present a systematic description for derivation of Algorithm~\ref{alg:async-admm}.
\subsection{Derivation of the proposed algorithm}\label{sec: alg derivation}
We can rewrite the distributed optimisation problem \eqref{eq: optim_problem}, as follows:
\begin{equation*}
    \begin{array}{ll}
         \underset{(\vx_1, \dots, \vx_n)}{\min} &\displaystyle{\sum_{i=1}^n f_i(\vx_i) + \sum_{j=1}^m \iota_{C_{\cV_j}}(\vx_{\cV_j}),}
    \end{array}
\end{equation*}
where $\vx_{\cV_i} := \textrm{col}(\{ \vx_j\}_{j\in \cV_i}) \in \RR^{d|\mathcal{V}_i|}$, for all $i\in [m]$ and $C_{\cV_i}:= \{\vx_{\cV_i} \mid \vx_\ell = \vx_k, \textrm{for all } (k, \ell) \in \cE_i\}$. By defining $\vz_i = \textrm{col}(\{ \vz_{i_j}\}_{j\in \mathcal{V}_i})$ we can rewrite again as below
\begin{equation} \label{eq: reformulated cp}
    \underset{(\vx, \vz)}{\min} \quad \displaystyle{\sum_{i=1}^n f_i(\vx_i) + g(\vz_1, \dots, \vz_m)}, \;     \text{subject to} \quad \vx_{\cV_j} = \vz_j.
\end{equation}
where $g(\vz_1, \dots, \vz_m) : = \sum_{j=1}^m \iota_{C_{\cV_j}}(\vz_j)$. Defining $\cT$ as a linear operator, $\cT(\vx_1, \dots, \vx_n) = (\vx_{\cV_1}, \vx_{\cV_2}, \dots, \vx_{\cV_m})$ simplifies the problem formulation as follows:
\begin{equation} \label{eq: reformulated cp subgraph}
\begin{array}{ccc}
    \underset{(\vx, \vz) \in \RR^d \times \RR^p}{\min} &\displaystyle{\sum_{i=1}^n f_i(\vx_i) + g(\vz)}, \\
    \text{subject to} & \vz = \cT(\vx).
\end{array}
\end{equation}
where $\vx= (\vx_1, \dots, \vx_n)$, $\vz = (\vz_1, \dots, \vz_m)$ and $p:= \sum_{j\in [m]} d|\cV_j|$. 

Revisiting the graph of Figure~\ref{fig: example graph}, we have $\cT(\vx) = (\vx_1, \vx_2, \vx_3, \vx_4, \vx_3, \vx_5, \vx_4, \vx_6)$ and sets $C_{\cV_i}, i=1,2, 3$ are defined as spaces of the form $C_{\cV_1} = \{(\va, \va, \va, \va) ; \va \in \RR^d\}$, $C_{\cV_2} = \{(\vb, \vb); \vb \in \RR^d\}, C_{\cV_3} = \{ (\vc, \vc); \vc \in \RR^d\}$. The function $g$ in this example is the indicator of vector space of the form $\{(\va,\va,\va,\va, \vb,\vb,\vc, \vc)\}$. Therefore, the function $g$ equals to zero if and only if $\vx_1=\vx_2=\vx_3=\vx_4, \vx_3 = \vx_5$ and $\vx_4=\vx_6$. Thus, all $\vx_i$'s are equal.

Now the problem \eqref{eq: reformulated cp subgraph} is of the form that we can apply ADMM which amounts to the following iterates
\begin{subequations}
    \begin{align}
        \vx^+ &\in \argmin \cL_\beta(\vx, \vy, \vz), \\
        \vz^+ &\in \argmin \cL_\beta(\vx^+, \vy, \vz),\\
        \vy^+ &= \vy + \beta(\cT(\vx^+) - \vz^+),
    \end{align}
\end{subequations}
where $\cL_\beta(\vx, \vz, \vy):= \sum_{i=1}^n f_i(\vx_i) + g(\vz) + \langle \vy, \cT (\vx) - \vz\rangle + \frac{\beta}{2} \|\cT(\vx) - \vz \|^2$, $\vy_i = \textrm{col}(\{ \vy_{i_j}\}_{j\in \mathcal{V}_i})$, and $\vy = (\vy_1, \dots, \vy_m)$. Define $\cN_i := \{j \mid i \in  \cV_j\}$. The update equations above can be written as
\begin{subequations} \label{eq: sync_admm_graph}
    \begin{align} 
        \vx_i^+ &\in \argmin f_i(\vx_i)\! +\sum_{j\in \cN_i} \langle \vy_{j_i} , \vx_i \rangle+ \frac{\beta}{2} \sum_{j\in\cN_i} \|\vx_i - \vz_{j_i} \|^2 \nonumber \\
         &= \argmin f_i(\vx_i) + \frac{\beta}{2} \sum_{j\in\cN_i}\|\vx_i - (\vz_{j_i} - \frac{1}{\beta}\vy_{j_i})\|^2 \nonumber \\
         & = \prox_{\frac{1}{\beta |\cN_i|} f_i} \left(\frac{1}{|\cN_i|}\sum_{j\in\cN_i} (\vz_{j_i} - \frac{1}{\beta}\vy_{j_i})\right), \: \forall i\in [n]\\
         \vz_i^+\! &\in \argmin \iota_{C_{\cV_i}}(\vz_i) + \frac{\beta}{2}\|\vz_i -((\cT\vx^+)_i +\frac{\vy_i}{\beta}) \|^2 \nonumber \\
          &\stackrel{(a)}{=}\! \mathbf{1}_{|\cV_i|} \otimes \overbrace{\frac{1}{|\cV_i|} \sum_{j\in \cV_i} (\vx_j^+ + \frac{\vy_{i_j}}{\beta})}^{=:\overline{\vz}_i^+}, \quad \forall i\in [m]\\
         \vy_{i_j}^+ &= \vy_{i_j} + \beta(\vx_j^+ - \bar{\vz}_i^+), \quad \forall i\in[m],\  \forall j\in \cV_i.
    \end{align}
\end{subequations}
where $(\cT(\vx^+))_i :=\vx^+_{\cV_i}$ and $\boldsymbol{1}_{|\cV_i|} \otimes \bar{\vz}_i^+ = (\bar{\vz}_i^+, \dots, \bar{\vz}_i^+)$ where vector $\bar{\vz}_i^+$ has been repeated $|\cV_i|$ times. and (a) is due to the fact that the proximal operator of the indicator function is the projection operator. Performing the equations above asynchronously in the sense that at each iteration $k$ one subgraph $i^k\in [m]$ or a set of subgraphs $\cS^k\subseteq [m]$ are chosen and all nodes in the subgraph perform the updates, results in Algorithm \ref{alg:async-admm}. 

In the following, we establish the almost sure convergence properties of the algorithm and provide a sufficient condition on the cardinality of $\cS^k$ (the set of nodes that update at each iteration) so that the algorithm converges.
\section{Convergence Analysis} \label{sec: proof BC-DRS}
\noindent For proving the convergence properties of Algorithm \ref{alg:async-admm} we first show that the algorithm is equivalent to the Randomised Block Coordinate DRS (RBC-DRS) method and then prove the convergence of RBC-DRS which implies the convergence of the proposed method. 
\subsection{Equivalence of Asynchronous ADMM and RBC-DRS}
Operator splitting methods in optimisation address the problems of the form
\begin{align} \label{eq: DRS problem}
    \begin{array}{ll}
    \underset{\vs \in \RR^p}{\min} &\displaystyle{\varphi(\vs):=\widetilde{f}(\vs) + \widetilde{g}(\vs)}.
    \end{array}
\end{align}

DRS method is one of the important splitting based algorithm in the literature. One step of DRS iteration with step-size $\gamma$ is as follows
\begin{subequations}
\begin{align}
    &\vu \in \prox_{\gamma \widetilde{f}}(\vs),\\
    &\vv \in \prox_{\gamma \widetilde{g}}(2\vu - \vs),\\
    &\vs^+ = \vs + (\vv - \vu).
\end{align}
\end{subequations}
Where the proximal operator $\prox_{\gamma f}$ is defined as $\prox_{\gamma f}(\vs):=\argmin_{x} f(\vx) + \frac{1}{2\gamma}\|\vx - \vs \|^2$. The ADMM algorithm, which is closely related to DRS, applies to linearly constrained problems of the form
\begin{equation} \label{eq: ADMM-contrained}
    \underset{(\vx, \vz)}{\min}  \quad \displaystyle{f(\vx) + g(\vz)}, \qquad \text{subject to} \quad \vz = M\vx.
\end{equation}

The steps of ADMM are outlined below:
\begin{subequations}
\label{ADMM-alg-2blocks}
\begin{align}
    &\vx^{+} \in \argmin \cL_{\beta}(\ \scalebox{0.6}{$\bullet$}\ , \vz, \vlambda),\\
    &\vlambda^+ = \vlambda + \beta (M\vx^{+} - \vz ) \label{eq: admm-lambda},\\
    &\vz^+ \in \argmin \cL_{\beta}(\vx^{+}, \ \scalebox{0.6}{$\bullet$}\ , \vlambda^{+}).    
\end{align}
\end{subequations}
where $\cL_{\beta}(\vx, \vz, \vlambda) := 
    f(\vx) + g(\vz) + \langle \vlambda, M\vx -\vz\rangle 
    +\frac{\beta}{2}\|M\vx - \vz  \|^2$
and the scalar $\beta$ is the \emph{penalty parameter}.

We are interested in establishing a relation between Randomised Block Coordinate DRS (RBC-DRS) and a version of ADMM. In the block coordinate DRS applied to \eqref{eq: DRS problem}, at each iteration one chooses block coordinate $i$ uniformly at random from $[m]$ and updates $i$-th block component of $\vs$ as follows
\begin{subequations}\label{RBC-DRS}
    \begin{align}
        &\vu \in \prox_{\gamma \widetilde{f}} (\vs), \\
        &\vv \in \prox_{\gamma \widetilde{g}} (2\vu - \vs), \\
        &\vs^+ = \vs + (\vv - \vu)_i\ve_i.
    \end{align}
\end{subequations} 
where $\ve_i$ is the vector that is one in $i$-th block component and zero elsewhere. Theorem below formalises the relation between Asynchronous ADMM and BC-DRS.

\begin{theorem} \label{thm: asynaddm-rbcdrs equiv}
    The Randomised Block-Coordinate Douglas-Rachford Splitting (RBC-DRS) method applied to \eqref{eq: DRS problem} is equivalent to an asynchronous decentralised ADMM algorithm given in Algorithm \ref{alg:async-admm} with $\vs :=M\vx-\vlambda/\beta,\vu:=M\vx, \vv:=\vz,$
        and $\widetilde{f} (\vs) := \inf_{\vx \in \{\vx \mid M\vx = \vs \}} f(\vx) , \widetilde{g} \equiv g$ and $\gamma = \frac{1}{\beta}$.
    
\end{theorem}

\begin{proof}
    At each iteration, a subset of nodes, $\cV_i$, is selected and $\vs^+ = \vs + (\vv - \vu)_{\cV_i}$. Let $\vs = (\vs_{\cV_1}, \dots, \vs_{\cV_m})$. Therefore,  $ \vs^+ = (\vs_{\cV_1}, \dots, \vs_{\cV_{i-1}}, \vs_{\cV_i}^+, \vs_{\cV_{i+1}}, \dots, \vs_{\cV_m})$. Let $\vs_{\cV_r}^+ = (\vs_{r_1}^+, \vs_{r_2}^+, \dots, \vs_{r_{|\cV_r|}}^+)$, $\vx_{\cV_r}^+ = (\vx_{r_1}^+, \vx_{r_2}^+, \dots, \vx_{r_{|\cV_r|}}^+)$, and $\vy_r^+ = (\vy_{r_1}^+, \dots, \vy_{r_{|\cV_r|}}^+)$. By the same change of variables in \cite[Thm ~5.5]{themelis2020douglas}, that is:
    \begin{subequations} \label{eq: change of variables}
    \begin{alignat}{3}
        &\vs :=M\vx-\vlambda/\beta, &&\vu:=M\vx, & &\vv:=\vz,\\
        &\vs^+:=M\vx^+-\vlambda^+/\beta,\  &&\vu^+:=M\vx^+,\  & &\vv^+:=\vz^+,
    \end{alignat}
    \end{subequations}

    we have $\vs_{\cV_j}^+ = \vs_{\cV_j} = \vx_{\cV_j} - \vy_j/\beta$ for $j \neq i$ and $\vs_{\cV_i}^+ =  \vx_{\cV_i}^+ - \vy_i^+/\beta$ where $\vx$ and $\vy$ are updated according to \eqref{eq: sync_admm_graph}. For $j\neq i$ and $k \in [|\cV_j|]$ one has $(\vs^+)_{j_k} = (\vs)_{j_k} = (M\vx - \vy/\beta)_{j_k} = (M\vx^+ - \vy^+/\beta)_{j_k}$. For components $\vs_{i_k}$ where $k\in [|\cV_i|]$ one has $(\vs^+)_{i_k} = (\vs)_{i_k} + (\vv - \vu)_{i_k}$. Then
\begin{align*}
    (\vs^+)_{i_k} &= (M\vx - \vy/\beta)_{i_k} + (\vz - M\vx)_{i_k} \\
    &\stackrel{(a)}{=} (M\vx - \vy/\beta)_{i_k} + (\vy/\beta - \vy^+/\beta +M\vx^+)_{i_k}\\
    & = (M\vx^+ - \vy^+/\beta)_{i_k},
\end{align*}
where in (a) we have used \eqref{eq: admm-lambda} which gives $\vz = \vy/\beta - \vy^+/\beta + M\vx^+$. This concludes the proof.
\end{proof}
\subsection{Convergence of Randomised Block Coordinate DRS method}
In this section, we analyse the block coordinate DRS method under the following standing and standard assumptions (see Section \ref{notations and preliminaries} for definitions of smooth, proper, lower semi-continuous, and prox-bounded functions):
\begin{assumption}\label{assp: optim problem}
The following hold:
    \begin{enumerate}
        \item The function $\varphi(\vs)$ is lower bounded and the  problem \eqref{eq: DRS problem} has a solution.
        \item The function $\tilde{f}$ is $L$-smooth and $\tilde{f}+ \frac{\ell}{2}\| \cdot \|^2$ is convex for some real number $\ell \in \RR$.
        \item The function $\tilde{g}$ is proper, lower semicontinuous, and prox-bounded.
    \end{enumerate}
\end{assumption}

By assuming that at each iteration of RBC-DRS block $\cC^k \subseteq [p]$ of the optimisation variable is updated, we need the following assumption to ensure each component is updated frequently enough. It is important to note that the uniformity of the distribution does not play an important role in the proof and one only needs $\min_{j\in [p]} \PP(j\in \cC^k) >0$, see Remark \ref{poe}. We consider the uniform distribution for the sake of the simplicity of the analysis.

\begin{assumption}\label{assp: Random sets} At each iteration $k$, the set
$\cC^k \subseteq [p]$ is chosen uniformly at random, i.e. $\PP(j\in \cC^k) =
\frac{1}{p}$ for all $j\in [p]$. 
% \footnote{The corresponding probability space, adopted from the model in \cite{cannelli2020asynchronous}, $(\Omega, \cF, \PP)$ is defined in Appendix \ref{appx: prob model}.}.
\end{assumption}

The main result is presented below.
\begin{theorem} \label{thm: main result}
    Let the functions $\varphi$, $\tilde{f}$ and $\tilde{g}$ satisfy Assumption \ref{assp: optim problem}, the index sampling $\cC^k$ follows the assumption \ref{assp: Random sets} in the Douglas-Rachford Splitting method \eqref{RBC-DRS}. If $|\cC^k| = |\cC|$ for all $k$ where $\frac{p}{2} + \delta < |\cC| \leq p$ for some $\delta>0$ and $\gamma>0$ is sufficiently small such that
    $$(1+\gamma L)^2 + \gamma L^2 + \frac{5\gamma \ell}{2} + (1+\gamma)(1+\gamma \ell)^2(1-\frac{|\cC|}{p})-\frac{3}{2} < 0, $$
    then the following hold:
    \begin{enumerate}[label=\roman*.]
        \item The limit of the random sequence $\cL^k:= \cL_{\gamma}(\vu^k, \vv^k, \vs^k)$ exists and is bounded with probability one, i.e., $\lim_{k\to \infty} \cL^k \stackrel{(a.s)}{=} \cL^* < \infty$.
        \item The sequence of random variables $\|\vv^k - \vu^k\|^2$ is summable with probability one.% (it vanishes in particular).
        \item The random sequences $\vv^k$ and $\vu^k$ have the same almost sure limit points and all of which are stationary points of the objective function $\varphi$.
    \end{enumerate}
    where $\cL_\gamma(\vu, \vv, \vs):= \tilde{f}(\vu) + \tilde{g}(\vv)   +\frac{1}{\gamma} \langle \vs - \vu, \vv - \vu\rangle -\frac{1}{2\gamma} \|\vv - \vu \|^2$.
\end{theorem}
Before proving the theorem, we present some useful lemmas.
We consider the general inexact iterations for the DRS method where at each iteration there is some source of noise in $\vs$-update step. Two consecutive iterations of the inexact DRS algorithm are as follows.
\begin{equation}\label{eq: DRS algorithm-1}
\begin{aligned}
        \vu &\in \prox_{\gamma \tilde{f}} (\vs^-),         && \vu^+ \in \prox_{\gamma \tilde{f}} (\vs), \\
        \vv &\in \prox_{\gamma \tilde{g}} (2\vu - \vs^-),  && \vv^+ \in \prox_{\gamma \tilde{g}} (2\vu^+ - \vs), \\
        \vs &= \vs^- + (\vv - \vu) + \vxi,     &&\vs^+ = \vs + (\vv^+ - \vu^+) + \vxi^+,
\end{aligned}
\end{equation}
where random sequence $\vxi$ captures the effect of randomness caused by the block coordinate update of the algorithm in each iteration, namely $\vxi_k = \frac{p}{|\cC^k|}\sum_{j\in \cC^k} [\vv^k - \vu^k]_j\ve_j - (\vv^k - \vu^k)$ where $\cC^k \subseteq [p]$ and $[\va]_j$ denotes the $j$-th block component of vector $\va$. 

We use \textit{Douglas-Rachford merit function} defined as
\begin{align} \label{eq: DR merit function}
    \cL_{\gamma}(\vu, \vv, \vs)&:=\! \tilde{f}(\vu) + \tilde{g}(\vv)   +\frac{1}{\gamma} \langle \vs - \vu, \vv - \vu\rangle -\frac{1}{2\gamma} \|\vv - \vu \|^2\\
    & \stackrel{(a)}{=}\! \tilde{f}(\vu) + \tilde{g}(\vv) +\frac{1}{2\gamma} \|2\vu - \vs - \vv \|^2  \nonumber \\
    &\qquad - \frac{1}{2\gamma}\|\vs - \vu \|^2 
    -\frac{1}{\gamma} \|\vv - \vu \|^2 \label{eq: DR merit 2} \\
    & \stackrel{(b)}{=}\! \tilde{f}(\vu)\! +\! \tilde{g}(\vv)\! +\!\frac{1}{2\gamma}\!\left(\|\!\vs\! -\! \vu\! \|^2\! -\! \|\! \vs\! -\! \vv\! \|^2\! \right), \label{eq: DR merit 3}
\end{align}
where $\langle\! \va\!,\! \vb\! \rangle = \frac{1}{2} (\|\va\! +\! \vb\|^2\! - \|\va\|^2\! -\! \|\vb\|^2)$ and  $\langle\! \va\!,\! \vb\!  \rangle = \frac{1}{2} (\!\|\va\|^2\! +\! \|\vb\|^2\! - \|\va\! -\! \vb\|^2)$ are used in (a) and (b), respectively.

The following lemmas are needed for establishing the result.
\begin{lemma}\label{lemma: differences}
    Let $\tilde{f}$ and $\tilde{g}$ satisfy Assumption \ref{assp: optim problem}, $\cL_{\gamma}$ be defined as \eqref{eq: DR merit function}, and $\vu, \vv, \vs$ follow \eqref{eq: DRS algorithm-1}. Then,
\begin{enumerate}[wide, labelindent=0pt, label=\roman*.]
        \item  $\cL_\gamma(\vu^+, \vv^+, \vs^+)\! -\cL_{\gamma}(\vu^+, \vv^+, \vs)\! =\! \frac{1}{\gamma}\langle \vs^+ - \vs, \vv^+ - \vu^+\rangle$.
        \item $\cL_\gamma(\vu^+, \vv^+, \vs) -\cL_{\gamma}(\vu^+, \vv, \vs) \leq -\frac{1}{\gamma} \|\vv^+ - \vu^+ \|^2 + \frac{1}{\gamma}\| \vv - \vu \|^2+ \frac{1}{\gamma}\| \vu^+ - \vu \|^2 -\frac{2}{\gamma}\langle \vv - \vu, \vu^+ - \vu \rangle $.
        \item $\cL_\gamma(\vu^+, \vv, \vs) -\cL_{\gamma}(\vu^+, \vv, \vs) \leq -\frac{1}{2}(\frac{1}{\gamma} - \ell)\|\vu^+ - \vu \|^2$.
    \end{enumerate}
\end{lemma}

\begin{proof}

    (i) follows from the definition of the Douglas-Rachford merit function in \eqref{eq: DR merit function}. Using \eqref{eq: DR merit 2} one has
    \begin{align*}
        &\cL(\vu^+, \vv^+, \vs) - \cL(\vu^+, \vv, \vs) \\
        & = \left(\tilde{g}(\vv^+) + \frac{1}{2\gamma}\|2\vu^+ - \vs - \vv^+ \|^2 - \frac{1}{\gamma} \|\vv^+ - \vu^+ \|^2\right) \\
        &\quad - \left(\tilde{g}(\vv) + \frac{1}{2\gamma}\|2\vu^+ - \vs - \vv \|^2 - \frac{1}{\gamma} \|\vv - \vu^+ \|^2\right)
    \end{align*}
    and using the fact that $\vv^+ \in \prox_{\gamma \tilde{g}} (2\vu^+ - \vs)$ we deduce that $\tilde{g}(\vv^+) + \frac{1}{2\gamma}\|2\vu^+ - \vs - \vv^+ \|^2 \leq \tilde{g}(\vv) + \frac{1}{2\gamma}\|2\vu^+ - \vs - \vv \|^2$ which implies that 
    \begin{align}
        &\cL(\vu^+, \vv^+, \vs) - \cL(\vu^+, \vv, \vs) \nonumber\\
        &\qquad\leq -\frac{1}{\gamma} \|\vv^+ - \vu^+ \|^2 + \frac{1}{\gamma} \|\vv - \vu^+ \|^2
    \end{align}
    using $\vv - \vu^+ = \vv - \vu - (\vu^+ - \vu)$ completes the proof of part (ii). Using \eqref{eq: DR merit 3} one has
    \begin{align*}
        &\cL(\vu^+, \vv, \vs) - \cL(\vu, \vv, \vs) \\
        &\quad= \left(\tilde{f}(\vu^+) + \frac{1}{2\gamma}\| \vs - \vu^+ \|^2 \right) - \left( \tilde{f}(\vu) + \frac{1}{2\gamma} \|\vs - \vu \|^2 \right).
    \end{align*}
    Using the convexity of $\tilde{f}(\cdot) + \frac{\ell}{2} \| \cdot \|^2$ we infer that $\tilde{f}(\cdot) + \frac{1}{2\gamma}\|\cdot - s\|^2$ is strongly convex with parameter $\frac{1}{\gamma} - \ell$ which in turn with the differentiability of $\tilde{f}$ implies that
    \begin{align}
        &\cL(\vu^+, \vv, \vs) - \cL(\vu, \vv, \vs) \nonumber\\
        &\quad \leq - \langle \nabla \tilde{f}(\vu^+) + \frac{1}{\gamma}(\vu^+ - \vs), \vu - \vu^+ \rangle \nonumber \\
        &\qquad -\frac{1}{2}(\frac{1}{\gamma} - \ell)\|\vu^+ - \vu \|^2.
    \end{align}
    Invoking the optimality condition of $\vu^+ \in \prox_{\gamma \tilde{f}}(\vs)$, that is $\vs = \vu^+ + \gamma \nabla \tilde{f}(\vu^+)$, proves part (iii).
\end{proof}
\begin{lemma} \label{lemma: L difference}
If the condition of Lemma \ref{lemma: differences} are satisfied then 
\begin{align*}
    &\cL(\vu^+ , \vv^+, \vs^+)-\cL(\vu, \vv, \vs) 
    \\
    & \leq \frac{1}{\gamma}\langle \vs^+ - \vs, \vv^+ - \vu^+\rangle -\frac{1}{\gamma} \|\vv^+ - \vu^+ \|^2 + \frac{1}{\gamma}\| \vv - \vu \|^2 \\
    & + \frac{1}{\gamma}\| \vu^+\!\! -\! \vu \|^2 \!
    -\!\frac{2}{\gamma}\langle\! \vv\! -\! \vu, \vu^+\!\! -\! \vu \rangle \!-\!(\frac{1}{2\gamma} - \frac{\ell}{2})\|\vu^+\!\! -\! \vu \|^2.
\end{align*}
\end{lemma}
\begin{proof}
    Adding up the results of parts (i), (ii), and (iii) of Lemma \ref{lemma: differences} proves the result.
\end{proof}
\begin{lemma}\label{lemma: xi var bound}
    Let $\cC^k$ be a  subset of $[p]$ chosen uniformly at random, then $        \EE(\|\vxi_k \|^2 \mid \cF_{k-1}) \stackrel{(a.s.)}{\leq} \zeta \EE( \|\vu^{k+1} - \vu^k \|^2\mid \cF_{k-1})$,
where $\zeta := (1 - \frac{|\cC^k|}{p})(1+\gamma L)^2$ and $\cF_k$ is the $\sigma$—algebra generated by randomness up to time $k$.
\end{lemma}
\begin{proof}
\begin{align*}
    &\EE(\|\vxi_k \|^2 \mid \cF_{k-1}) \\
    &\quad  = \EE \left\{ \| \frac{p}{|\cC^k|}\sum_{j\in \cC^k} [\vv^k-\vu^k]_j\ve_j - (\vv^k - \vu^k) \|^2 \mid \cF_{k-1} \right\} \\
    & \quad = \EE \left\{ \frac{p^2}{|\cC^k|^2}\sum_{j\in \cC^k} |[\vv^k-\vu^k]_j|^2 + \|\vv^k - \vu^k \|^2 \right. \\
    &\left.\qquad-\frac{2p}{|\cC^k|}\sum_{j\in |\cC^k|}|[\vv^k-\vu^k]_j|^2 \mid \cF_{k-1}\right\} \\
    & \quad \stackrel{(a)}{=} (\frac{p^2}{|\cC^k|^2} -\frac{2p}{|\cC^k|})\EE(\sum_{j\in \cC^k} |[\vv^k-\vu^k]_j|^2 \mid \cF_{k-1}) \\
    &\qquad+ \|\vv^k - \vu^k \|^2 \\
    & \quad= \left(\frac{|\cC^k|}{p}(\frac{p^2}{|\cC^k|^2} -\frac{2p}{|\cC^k|}) + 1\right)\|\vv^k - \vu^k \|^2 \\
    &\quad = (\frac{p}{\cC^k}-1)\|\vv^k - \vu^k \|^2 \\
    & \quad= (\frac{p}{|\cC^k|}-1)\frac{|\cC^k|}{p}\EE(\|\vs^k - \vs^{k-1} \|\mid \cF_{k-1}) \\
    & \quad= (1 - \frac{|\cC^k|}{p})\EE(\|\vs^k - \vs^{k-1} \|\mid \cF_{k-1}) \\
    & \quad \stackrel{(c)}{\leq} (1 - \frac{|\cC^k|}{p})(1+\gamma L)^2 \EE( \|\vu^{k+1} - \vu^k \|^2\mid \cF_{k-1}),
\end{align*}
    where in (a) we have used the face that $\vu^k, \vv^k$ are $\cF_k$-measurable random variables, (b) follows from taking conditional expectation with respect to $\cF_{k-1}$ from both sides of $\vs^{k} - \vs^{k-1} = \frac{p}{|\cC^k|} \sum_{j\in \cC^k}[\vv^k - \vu^k]_j\ve_j$ and in (c) we have used part (2) of Lemma \ref{lemma: prox coercivity} (see Appendix \ref{sec: background}) , Cauchy-Schwartz inequality, and monotonicity of conditional expectation.
\end{proof}

\emph{Proof of Theorem \ref{thm: main result}}: rewriting the result of Lemma \ref{lemma: L difference} for iteration counter $k$ gives
\begin{align*}
    \cL^{k+1}& - \cL^{k} \\
    &\leq \frac{1}{\gamma}\langle \vs^{k+1} - \vs^{k}, \vv^{k+1} - \vu^{k+1}\rangle -\frac{1}{\gamma} \|\vv^{k+1} - \vu^{k+1} \|^2 \\
    & + \frac{1}{\gamma}\| \vv^k - \vu^k \|^2 + \frac{1}{\gamma}\| \vu^{k+1} - \vu^k \|^2 \\
    &-\frac{2}{\gamma}\langle \vv^k - \vu^k, \vu^{k+1} - \vu^k \rangle-\frac{1}{2}(\frac{1}{\gamma} - \ell)\|\vu^{k+1} - \vu^k \|^2,
\end{align*}
where $\cL^k := \cL(\vu^k, \vv^k, \vs^k)$. Considering the iteration $\vs^k = \vs^{k-1} + (\vv^k - \vu^k) + \vxi_k$, let $\cF_k := \sigma(\vxi_1, \vxi_2, \dots, \vxi_k)$. We have
\begin{align*}
    &\cL^{k+1} - \cL^{k} \\
    &\leq \frac{1}{\gamma} \langle \vv^{k+1} - \vu^{k+1} + \vxi_{k+1}, \vv^{k+1} - \vu^{k+1}\rangle\\
    &\quad -\frac{1}{\gamma} \|\vv^{k+1} - \vu^{k+1} \|^2 + \frac{1}{\gamma}\| \vv^k - \vu^k \|^2 
    + \frac{1}{\gamma}\| \vu^{k+1} - \vu^k \|^2 \\
    &\quad -\frac{2}{\gamma}\langle \vs^k - \vs^{k-1} - \vxi_k, \vu^{k+1} - \vu^k \rangle \\
    &\qquad-\frac{1}{2}(\frac{1}{\gamma} - \ell)\|\vu^{k+1} - \vu^k \|^2 \\
    & = \frac{1}{\gamma} \langle \vxi_{k+1}, \vv^{k+1} - \vu^{k+1}\rangle + \frac{1}{\gamma}\| \vv^k - \vu^k \|^2  \\
    &\qquad-\frac{2}{\gamma}\langle \vs^k - \vs^{k-1} - \vxi_k, \vu^{k+1} - \vu^k \rangle \\
    &\qquad +\frac{1}{2}(\frac{1}{\gamma} + \ell)\|\vu^{k+1} - \vu^k \|^2.
\end{align*}

Note that the inequality above holds surely (not almost surely).
    
We now take the conditional expectation of the terms on both sides of the inequality above and using the monotonicity of conditional expectation and the fact that $\vu^{k+1}$ and $\vv^{k+1}$ are $\cF_k$ measurable random variables we have
\begin{align*}
    &\EE(\cL^{k+1} - \cL^k \mid \cF_k ) \\
    &\stackrel{(a.s.)}{\leq} \frac{1}{\gamma}\| \vv^k - \vu^k \|^2 -\frac{2}{\gamma}\langle \vs^k - \vs^{k-1} - \vxi_k, \vu^{k+1} - \vu^k \rangle \\
    &\qquad+\frac{1}{2}(\frac{1}{\gamma} + \ell)\|\vu^{k+1} - \vu^k \|^2.
\end{align*}
Invoking part (i) of Lemma \ref{lemma: prox coercivity} and the equation $\vs^k = \vs^{k-1} + (\vv^k - \vu^k) + \vxi_k$ result in 
\begin{align*}
    &\EE(\cL^{k+1} - \cL^k \mid \cF_k ) \stackrel{(a.s.)}{\leq} \frac{1}{\gamma}\| \vv^k - \vu^k \|^2 +\frac{2}{\gamma}\langle  \vxi_k, \vu^{k+1} - \vu^k \rangle \\
    &\qquad +\left[\frac{2}{\gamma}(\gamma \ell-1) +\frac{1}{2}(\frac{1}{\gamma} + \ell) \right]\|\vu^{k+1} - \vu^k \|^2 \\
    & \leq \frac{1}{\gamma}\|\vs^k - \vs^{k-1} - \vxi_k \|^2 +\frac{2}{\gamma}\langle  \vxi_k, \vu^{k+1} - \vu^k \rangle \\
    & \qquad + \frac{1}{\gamma}(\frac{5\gamma \ell - 3}{2})\|\vu^{k+1} - \vu^k \|^2 \\
    & = \frac{1}{\gamma}\|\vs^k - \vs^{k-1}  \|^2 + \frac{1}{\gamma}\|\vxi_k\|^2 -\frac{2}{\gamma} \langle  \vxi_k, \vs^k - \vs^{k-1} \rangle \\
    &\qquad +\frac{2}{\gamma}\langle  \vxi_k, \vu^{k+1} - \vu^k \rangle + \frac{1}{\gamma}(\frac{5\gamma \ell - 3}{2})\|\vu^{k+1} - \vu^k \|^2 \\
    & \stackrel{(a)}{\leq} \frac{1}{\gamma}\left((1+\gamma L)^2 + \frac{5\gamma \ell - 3}{2}\right) \|\vu^{k+1} - \vu^k \|^2 \\
    &\qquad + \frac{2}{\gamma}\langle  \vxi_k, \vu^{k+1} - \vu^k - (\vs^k - \vs^{k-1}) \rangle + \frac{1}{\gamma}\|\vxi_k\|^2
\end{align*}
Where in (a) we have used the part (ii) of Lemma \ref{lemma: prox coercivity} (see Appendix \ref{sec: background}). Using Lemma \ref{lemma: prox properties} in the first iteration of the DRS algorithm yields
\begin{align*}
    &\EE(\cL^{k+1} - \cL^k \mid \cF_k ) \\
     &\qquad \stackrel{(a.s.)}{\leq} \frac{1}{\gamma}\left((1+\gamma L)^2 + \frac{5\gamma \ell - 3}{2}\right) \|\vu^{k+1} - \vu^k \|^2 \\
    &\qquad+ 2\langle  \vxi_k, \nabla \tilde{f}(\vu^{k+1}) - \nabla \tilde{f}(\vu^k) \rangle + \frac{1}{\gamma}\|\vxi_k\|^2
\end{align*}
Using Young's inequality, $\langle \va, \vb \rangle \leq\frac{1}{2} (\|\va\|^2 + \| \vb \|^2)$, and Lipschitz continuity of the gradient of $\tilde{f}$ amounts to
\begin{align*}
        &\EE(\cL^{k+1} - \cL^k \mid \cF_k ) \\
        &\quad \stackrel{(a.s.)}{\leq} \frac{1}{\gamma}\left((1+\gamma L)^2 + \gamma L^2+ \frac{5\gamma \ell - 3}{2}\right) \|\vu^{k+1} - \vu^k \|^2 \\
        &\qquad + (1 + \frac{1}{\gamma})\|\vxi_k\|^2. 
\end{align*}
 Next we take conditional expectation of both sides of the inequality above invoking the assumption \ref{assp: Random sets} that at iteration $k$, the algorithm involves updating a random subset of coordinates, $\cC^k\subseteq [p]$, and the definition of $\vxi_k = \frac{p}{|\cC^k|}\sum_{j\in \cC^k} [\vv^k - \vu^k]_j\ve_j - (\vv^k - \vu^k)$.
 Using monotonicity of conditional expectation and double expectation law (a.k.a the tower property), that is $\EE(\EE(X \mid \cG_2) \mid \cG_1) = \EE(X\mid \cG_1)$ for sigma algebras $\cG_1 \subseteq \cG_2$, we have
\begin{align} \label{eq: final bound for xi}
        &\EE(\cL^{k+1} - \cL^k \mid \cF_{k-1} ) \nonumber\\
        &\stackrel{(a.s.)}{\leq} \frac{1}{\gamma}((1+\gamma L)^2\! +\! \gamma L^2\! +\! \frac{5\gamma \ell\! -\! 3}{2}) \EE(\|\vu^{k+1}\! -\! \vu^k \|^2\!\mid\!\cF_{k-1}) \nonumber\\
        &\qquad+ (1 + \frac{1}{\gamma})\EE(\|\vxi_k\|^2 \mid \cF_{k-1}). 
\end{align}
% The following lemma deals with the last term in the RHS of the above inequality.

Using Lemma \ref{lemma: xi var bound} in inequality \eqref{eq: final bound for xi} yields
\begin{align}\label{mideq1}
    \EE(\cL^{k+1}  \mid \cF_{k-1}) \stackrel{(a.s.)}{\leq}  \cL^k + \alpha \EE( \|\vu^{k+1} - \vu^k \|^2\mid \cF_{k-1}),
\end{align}
where 
\begin{align*}
        &\alpha := \frac{1}{\gamma}\left[(1+\gamma L)^2 + \gamma L^2 +\frac{5\gamma \ell}{2} +\right. \\
        &\qquad \left. (1+\gamma)(1+\gamma \ell)^2 (1-\frac{|\cC^k|}{p})-\frac{3}{2} \right]. 
\end{align*}
If $|\cC^k| \equiv |\cC| >
\frac{p}{2} + \delta$ for some $\delta > 0$ and $\gamma$ is sufficiently small,
it can be easily seen that $\alpha$ is negative. By subtracting $\underline{\cL} = \inf_k{\cL^k}$ from both sides of \eqref{mideq1} and invoking Lemma \ref{lemma: super MG thm} (see Appendix \ref{sec: background}) results in the almost sure existence and finiteness of $\lim_{k\to
\infty} \cL_k$ which proves part (i). It also implies that $\sum_k
\EE\{\|\vu^{k+1} - \vu^{k} \|^2 \mid \cF_{k-1}\} \stackrel{(a.s)}{<} \infty$
which in turn gives $\lim_{k\to \infty} \EE (\| \vu^{k+1} - \vu^{k}\|^2 \mid
\cF_{k-1}) \stackrel{(a.s.)}{=} 0$. Using Markov inequality for every $\epsilon >0$ we have $\PP(\|\vu^{k+1} - \vu^k \| > \epsilon) \leq \frac{\EE \|\vu^{k+1} - \vu^k \|}{\epsilon^2}$. Summing over $k$ yields $\sum_{k=0}^\infty \PP (\|\vu^{k+1} - \vu^k \| > \epsilon) \leq \frac{\EE \sum_{k=0}^\infty \|\vu^{k+1} - \vu^k \|^2}{\epsilon^2} < \infty$.  Invoking Borel-Cantelli lemma \cite[Thm 1.27]{karr1993Probability} proves that $\|\vu^{k+1} - \vu^k \|$ converges to zero almost surely.
By using Lemma \ref{lemma: prox coercivity} (see Appendix \ref{sec: background}) 
we can see that $\sum_k \EE\{\|\vs^{k} - \vs^{k-1} \|^2 \mid \cF_{k-1}\} \stackrel{(a.s)}{<} \infty$ and thus $\lim_{k\to \infty} \EE (\| \vs^{k} - \vs^{k-1}\|^2 \mid \cF_{k-1}) \stackrel{(a.s.)}{=} 0$ and $\lim_{k\to \infty} \|\vs^{k} - \vs^{k-1} \| \stackrel{a.s.}{=}0$. Since $\vs^k - \vs^{k-1} = \frac{p}{|\cC^k|}\sum_{j\in \cC^k}[\vv^k - \vu^k]_j\ve_j$, we can conclude that the sequences $\vv^k$ and $\vu^k$ have the same cluster points, i.e. a subsequence of $(\vv^k)_{k \in K}$ converges to $\vv^*$ iff so does the subsequence $(\vu^k)_{k\in K}$ (almost surely), this proves part (ii). From the equation $\vs^{k-1} = \vu^k + \gamma \nabla \tilde{f}(\vu^k)$ and the continuity of $\nabla \tilde{f}$  we conclude that $\lim_{K \ni k\to \infty} \vs^{k-1} = \vu^* + \gamma \nabla \tilde{f}(\vu^*) =: \vs^*$. Therefore, from lemma
\ref{lemma: prox properties} (see Appendix \ref{sec: background})  we have $\vu^* = \prox_{\gamma \tilde{f}} (\vs^*)$. Combining the first and the second equation in DRS algorithm implies $\vv^k \in \prox_{\gamma \tilde{g}}(\vu^k -\gamma \nabla \tilde{f}(\vu^k))$. Passing the limit over the subsequence $K$ gives $\vu^* = \vv^* = \lim_{K \ni k \to \infty} \vv^k \in \lim_{K \ni k \to \infty} \prox_{\gamma \tilde{g}} (\vu^k - \gamma \nabla \tilde{f}(\vu^k)) $. Using Lemma \ref{lemma: osc of prox} (see Appendix \ref{sec: background})  implies that
\begin{align*}
    \vv^* &\in \limsup_{\vu^k \to \vu^*} \prox_{\gamma \tilde{g}}(\vu^k - \gamma \nabla \tilde{f}(\vu^k))\\
    & \qquad\subseteq \prox_{\gamma \tilde{g}}(\vu^* - \gamma \nabla \tilde{f}(\vu^*))
\end{align*}
which gives $\vu^* \in \prox_{\gamma \tilde{g}} (\vu^* -\gamma \tilde{f}(\vu^*))$. This implies $\boldsymbol{0} \in \partial (\tilde{f} + \tilde{g})(\vu^*)$ and proves that the cluster point $\vu^*$ is a stationary point of the optimisation problem and completes the proof.
\hfill \qed
\begin{remark}\label{poe}
    The result of Theorem \ref{thm: main result} still holds if we relax the uniformity in Assumption \ref{assp: Random sets}  as $\min_{j} \PP(j\in \cC^k) > 0$ for all $k$, i.e., the probability of activation of the components is nonuniform but bounded away from zero.
\end{remark}
\subsection{Convergence of Algorithm \ref{alg:async-admm}}
Using the result of the last two section we can analyse Algorithm \ref{alg:async-admm}. Note that the set $S^k$ in Algorithm \ref{alg:async-admm} and $\cC^k$ are related in that $\cS^k$ specifies the indices of the activated subgraphs in $[m]$ while $\cC^k$ corresponds to the indices of the components of the node variables within these activated subgraphs, represented in the vector $\cT(\vx) := (\vx_{\cV_1}, \vx_{\cV_2}, \dots, \vx_{\cV_m})$. Note that the cardinality of $\cT(\vx)$ is $p=\sum_{j\in [m]}d|\cV_j |$.

Combining the results of Theorems \ref{thm: asynaddm-rbcdrs equiv} and \ref{thm: main result} yields the following theorem.

\begin{theorem}
    In Algorithm \ref{alg:async-admm}, if the penalty parameter $\beta$ and the set of subgraphs $\cS^k$ are chosen such that 
    \begin{align*}
        \frac{(\beta\!+\! L)^2}{\beta^2} + \frac{L^2}{\beta} + \frac{5\ell}{2\beta} + \frac{(\beta\!+\!1)(\beta\!+\! \ell)^2}{\beta^3}(1\!-\!\frac{d|\cS|}{p})-\frac{3}{2} < 0, 
    \end{align*} 
    then every limit point of the sequence $\{\vx_k, \vy_k, \vz_k \}_{k=0}^{\infty}$ is a stationary point of the constrained optimisation problem \eqref{eq: reformulated cp subgraph}.
\end{theorem}
\begin{proof}
    Given the facts that $\gamma = \frac{1}{\beta}$ and $\cC^k$ is the corresponding components of the activated nodes $\cS^k$ in the vector $\cT(\vx)$ which implies $|\cC| = d|\cS|$, the proof follows from the equivalence of Algorithm \ref{alg:async-admm} and RBC-DRS and Theorem \ref{thm: main result}.
\end{proof}

\section{Numerical Experiments} \label{sec: numerical}
\noindent \textbf{Phase Retrieval.} We evaluate the effectiveness of Algorithm \ref{alg:async-admm} in addressing the decentralised phase retrieval problem, where $n$ detectors try to reconstruct an unknown signal from noisy measurements. These measurements capture the magnitude of linear transformations of the signal, and the problem is relevant in various domains such as wireless communications, physical sciences, and applications like passive envelope detectors, optical imaging, and X-ray diffractive imaging \cite{wang2020decentralized}. 
Centralised processing may be impractical due to concerns regarding privacy, computational overhead, and the distributed nature of data collection using sensors. In this context, each detector $i$ has gathered $m_i$ noisy measurements, denoted as $b_{ij}$ where $j\in\{1,\dots,m_i\}$. These measurements represent the squared norm of the inner product between the discrete signal $\vx_0 \in \CC^d$ and a certain direction $\vt_{ij} = \Re(\vt_{ij})+ \jmath\Im(\vt_{ij}) \in \CC^d$, where $\jmath = \sqrt{-1}$ and $\CC$ is the set of complex numbers. Specifically, the measurements are assumed to follow the form $b_{ij} = |\langle \vx_0, \vt_{ij} \rangle|^2 + \nu_{ij}$ for all $j\in \{1, \dots, m_i\}$, where $\nu_{ij}$ denotes the noise in the power measurement. The objective is to utilise $\langle \vx_0, \vt_{ij} \rangle$ and recover both the phase and magnitude of $\vx_0$. To achieve this, a commonly used formulation is based on minimising the least squares error between the observation $b_{ij}$ and the square modulus of the transformation, $\langle \vx, \vt_{ij} \rangle$. Thus, the objective function of agent $i$ can be formulated as:
\begin{equation} \label{eq: phase retrieval local cost}
    f_i(\vx) = \frac{1}{m_i}\sum_{j=1}^{m_i} (b_{ij} - |\langle \vx, \vt_{ij} \rangle|^2)^2.
\end{equation}
This, in turn, leads to a nonconvex problem of the form \eqref{eq: optim_problem}. 

In our numerical experiment, we set $n=15$, $d=32$, and $m_i=30$ for all $i\in[n]$. The directions $\vt_{ij}$ and noise $\nu_{ij}$ are independently drawn from  $\cN(\boldsymbol{0}_{2d}, \frac{1}{2}\vI_{2d})$ and $\cN(0, 0.01^2)$, respectively. The communication graph took the form of a ring graph. Therefore, we have $m=n$ and 
\begin{equation} \label{eq: groups_example}
    \cV_i = 
    \begin{cases}
        \{i, i+1 \}, & \text{if } i< n, \\
        \{1, n\} & \text{if } i = n.
    \end{cases}
\end{equation}

\begin{figure}[ht]
    \centering
    \subfloat[]{
    % \begin{subfigure}[b]{0.45\textwidth}
        \includegraphics[width=.4\textwidth]{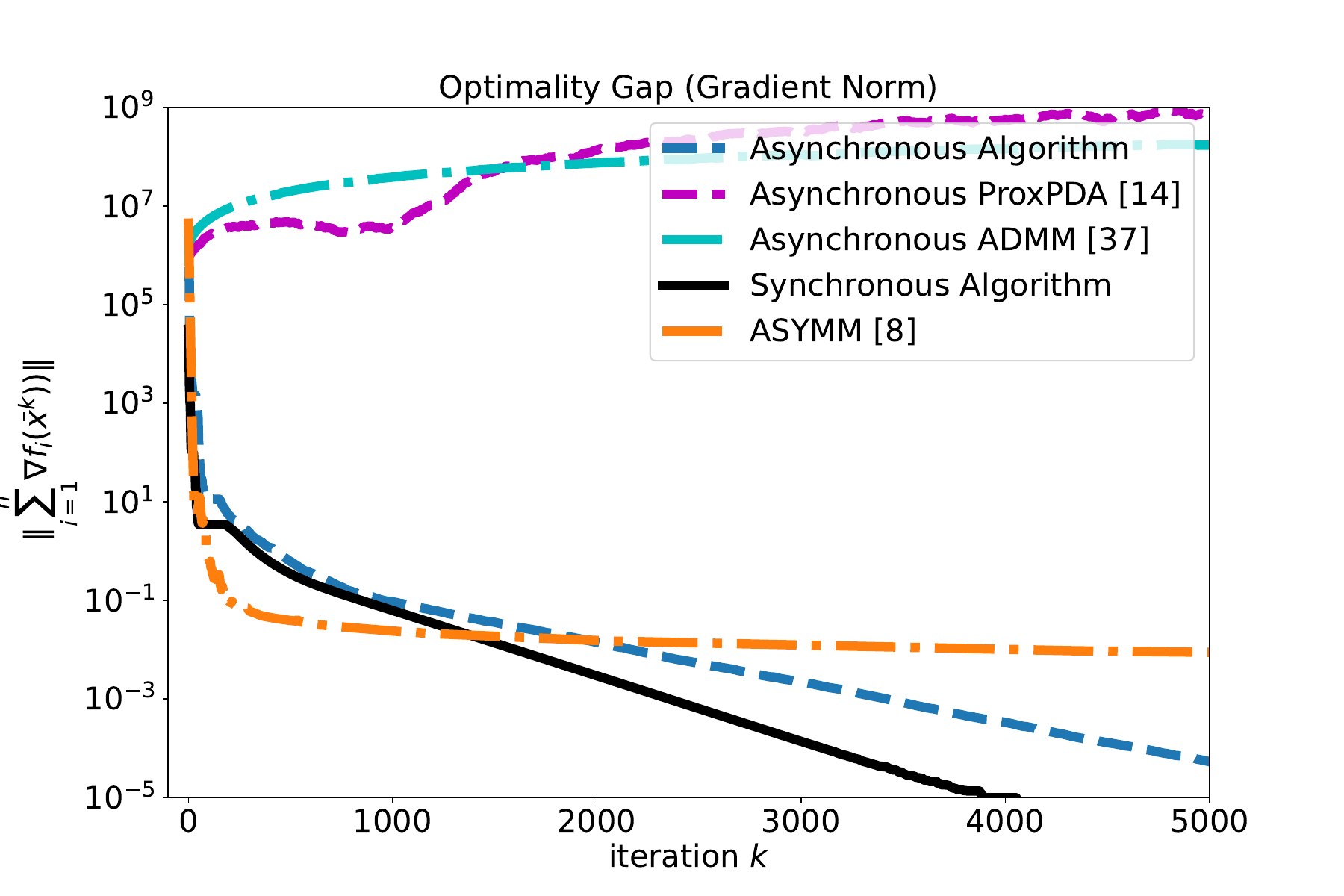}
        % \caption{$\log_{10}G^r(\bar{\vx}^r)$ vs. the number of iterations}
        \label{fig:D_r_20}}
    % \end{subfigure}
    \hfill
    % \begin{subfigure}[b]{0.45\textwidth}
    \subfloat[]{
        \includegraphics[width=.4\textwidth]{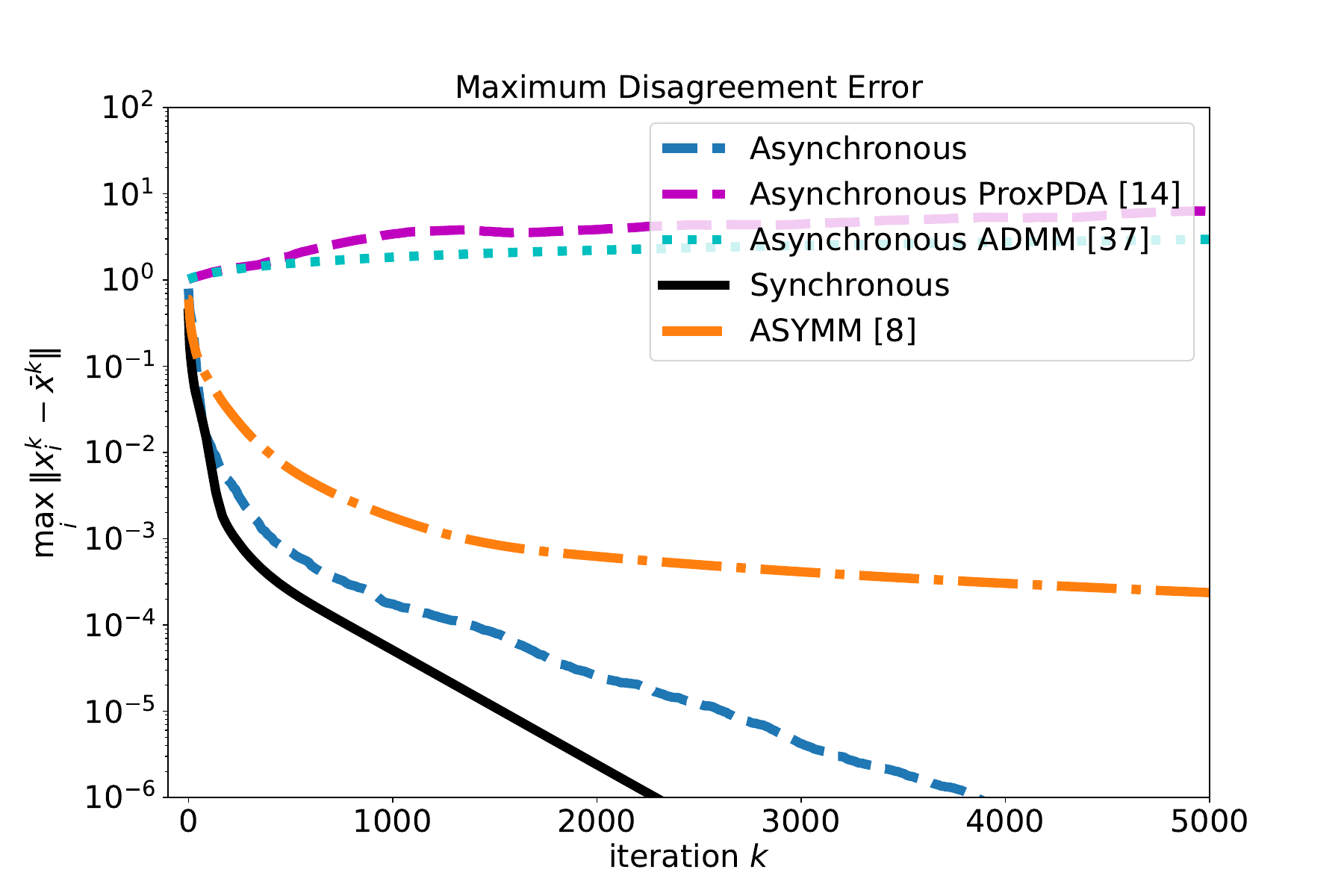}
        % \caption{$\log_{10}D^r(\vx^r)$ vs. the number of iterations}
        \label{fig:J_r_20}}
    % \end{subfigure}
    \caption{The performance of the proposed algorithm, ASYMM \cite{farina2019distributed}, asynchronous ProxPDA \cite{pmlr-v70-hong17a}, asynchronous distributed ADMM in \cite{yi2022sublinear}, and the synchronous version of the algorithm for solving the Phase Retrieval problem. (a) $\log_{10}G^r(\bar{\vx}^r)$ vs. the number of iterations, (b) $\log_{10}D^r(\vx^r)$ vs. the number of iterations.}
    \label{fig:PR vs ASYMM}
\end{figure}

To track the algorithms' convergence towards stationarity and consensus on the decision variable, we employ the gradient norm and disagreement gap. These are defined as $G^r:= \| \sum_{i=1}^n \nabla f_i(\bar{\vx}^r)\|$ and $D_r := \max_{i\in [n]} \|\vx_i^r - \bar{\vx}^r \|$, respectively. Furthermore, we have implemented the state-of-the-art Asynchronous Method of Multipliers (ASYMM) algorithm \cite{farina2019distributed} and asynchronous version (in the same sense as our proposed algorithm) of the algorithms in \cite{pmlr-v70-hong17a} and \cite{yi2022sublinear}. {From Figure \ref{fig:PR vs ASYMM} it is evident that both the synchronous and asynchronous proposed algorithms exhibit nearly identical performance, and both eventually outperform ASYMM and asynchronous versions of \cite{pmlr-v70-hong17a} and \cite{yi2022sublinear} diverge.  It should be noted that the last two algorithms are not designed for asynchronous settings, and their lack of convergence is therefore expected. We have conducted numerical simulations of these algorithms to highlight the nontrivial challenges associated with ensuring the convergence of asynchronous versions of distributed algorithms.} It is important to note that our definition of asynchrony differs from that in \cite{farina2019distributed}, where a more general type of asynchrony allows nodes to independently perform local updates without waking their neighbours.

\begin{figure}[ht]
    \centering
    \subfloat[]{
    % \begin{subfigure}[b]{0.45\textwidth}
        \includegraphics[width=.4\textwidth]{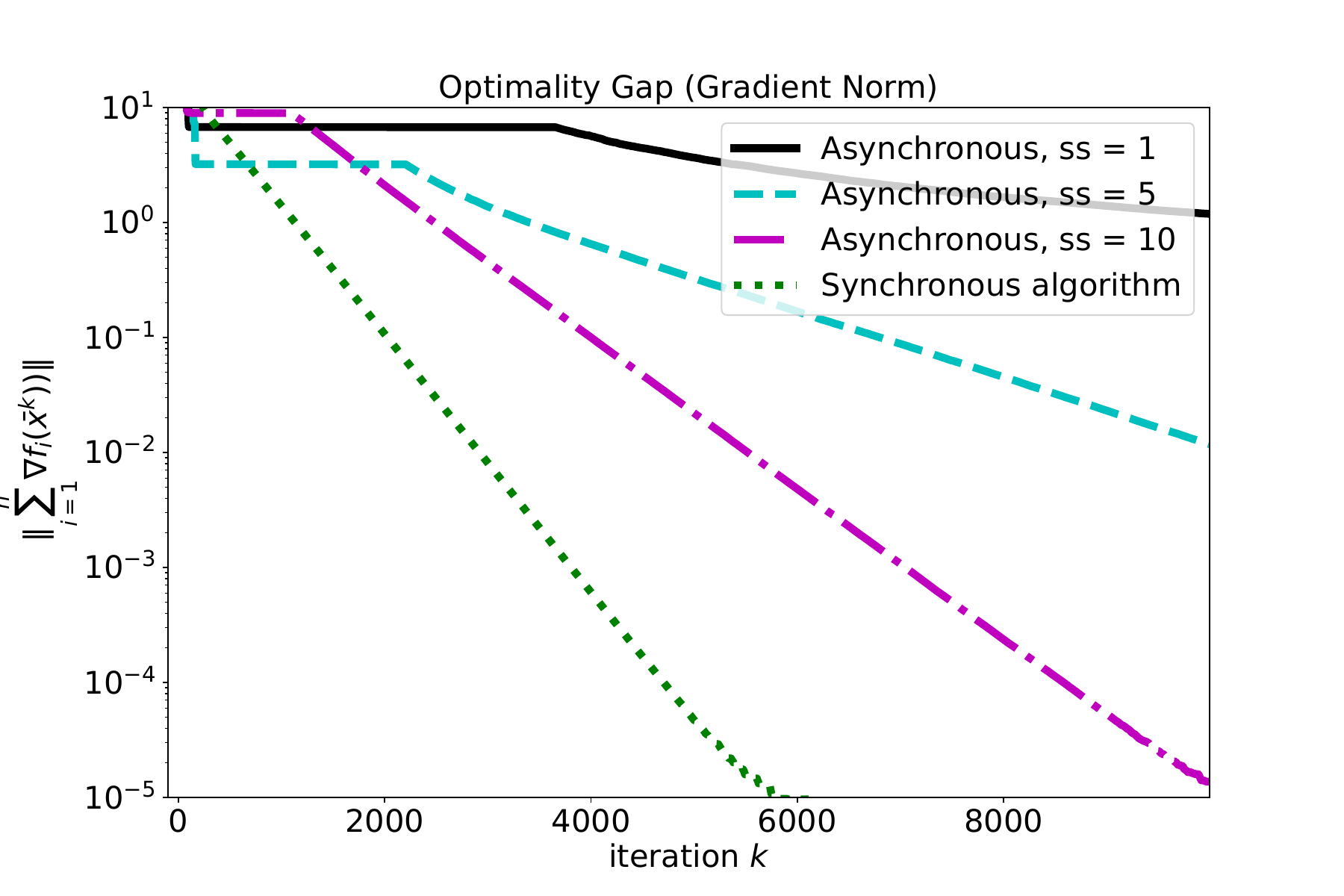}
        % \caption{$\log_{10}G^r(\bar{\vx}^r)$ vs. the number of iterations}
        \label{fig:G_PR_ss}}
    % \end{subfigure}
    \hfill
    % \begin{subfigure}[b]{0.45\textwidth}
    \subfloat[]{
        \includegraphics[width=.4\textwidth]{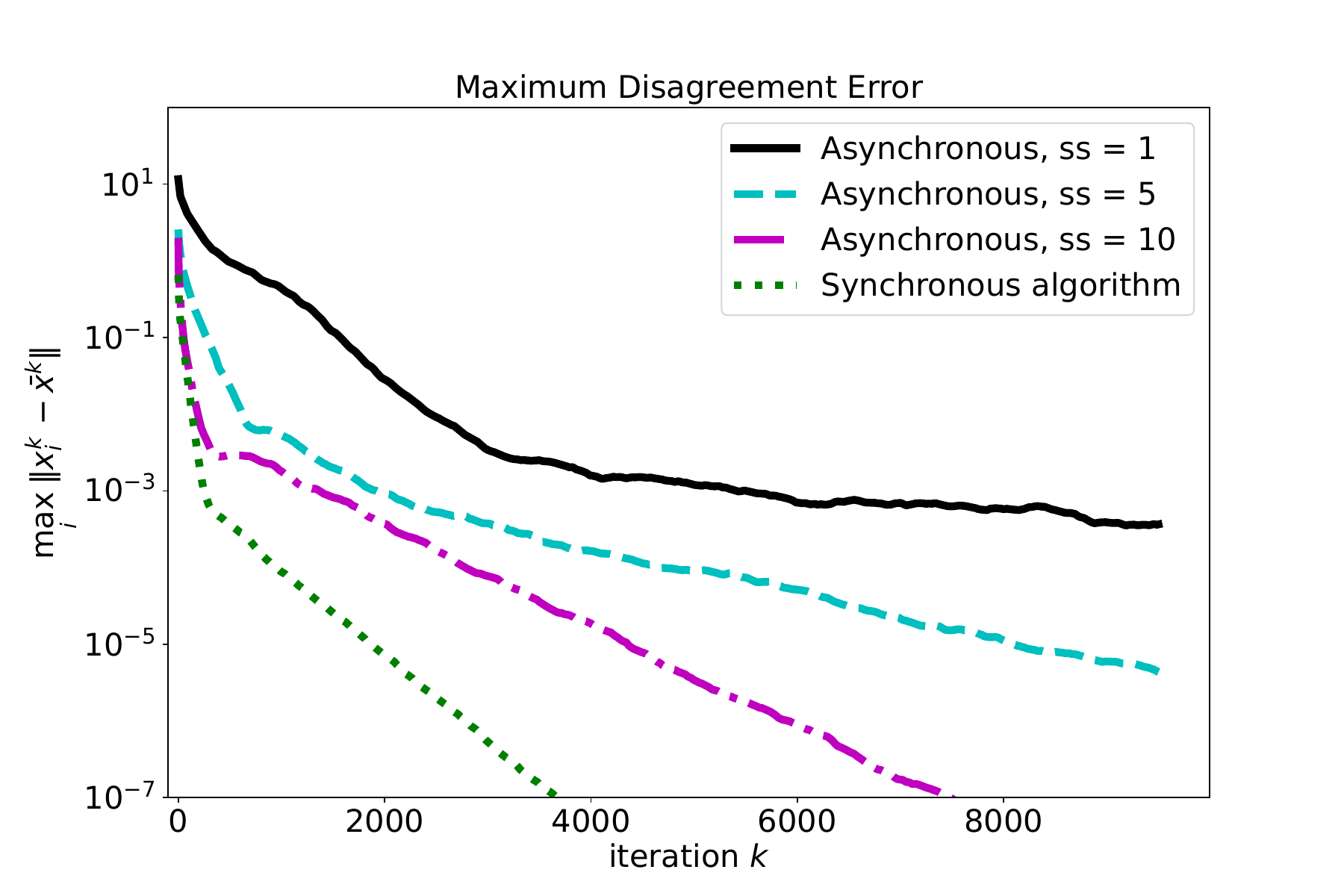}
        % \caption{$\log_{10}D^r(\vx^r)$ vs. the number of iterations}
        \label{fig:D_PR_ss}}
    % \end{subfigure}
    \caption{The performance of the proposed decentralised algorithm on the Phase Retrieval problem for different subset sizes. (a) $\log_{10}G^r(\bar{\vx}^r)$ vs. the number of iterations, (b) $\log_{10}D^r(\vx^r)$ vs. the number of iterations}
    \label{fig:PR}
\end{figure}

We examined the impact of the number of coordinates updated in each round of DRS iteration (or the number of agents activated in each round of asynchronous ADMM) in Figures \ref{fig:G_PR_ss} and \ref{fig:D_PR_ss} with $d=64$. The results indicate that as more nodes update at each iteration (subset size, denoted by $ss$ in the plots), the algorithm converges faster. Conversely, it is observed that insufficient activation of nodes in each iteration may lead to the algorithm not converging as desired.

\textbf{Sparse Principal Component Analysis.} Consider the sparse PCA problem \cite{richtarik2021alternating, chang2016asynchronous}:

\begin{equation} 
    \begin{array}{ccc}
                \underset{\vx \in \RR^d, \| \vx \|^2 \leq 1 }{\mbox{minimize}} &\displaystyle{\sum_{i=1}^n -\|P_i \vx \|^2 + \lambda \|\vx\|_1}, 
    \end{array}
\end{equation}
Here, $\lambda$ represents the regularization parameter for the $\ell_1$-penalised problem, and each agent $i$ independently possesses a data matrix $P_i \in \RR^{m_i \times d}$. Notably, existing algorithms in the literature, such as those outlined in \cite{hong2017distributed, chang2016asynchronous, scutari2019distributed}, are confined to scenarios with a central node in the communication graph, rendering them unsuitable for handling decentralised problems. Furthermore, each $f_i(\vx) :=-\vx^T P_i^T P_i \vx$ is a smooth concave function, facilitating a closed-form solution for the sub-problems when $\beta$ is sufficiently large.

\begin{figure}[ht]
    \centering
    \subfloat[]{
    % \begin{subfigure}[b]{0.45\textwidth}
        \includegraphics[width=.4\textwidth]{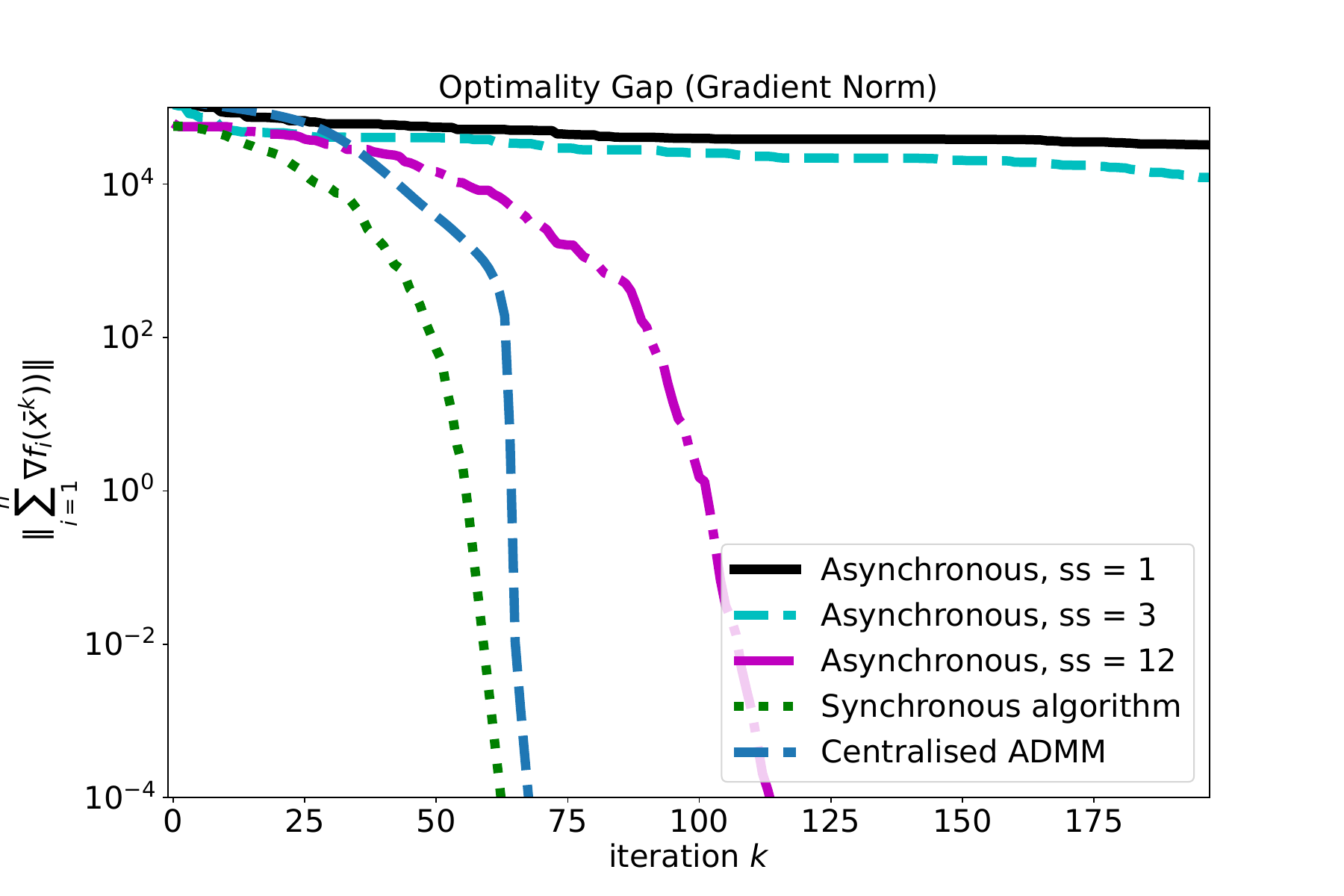}

        \label{fig:G_PCA}}
    % \end{subfigure}
    \hfill

    \subfloat[]{
        \includegraphics[width=.4\textwidth]{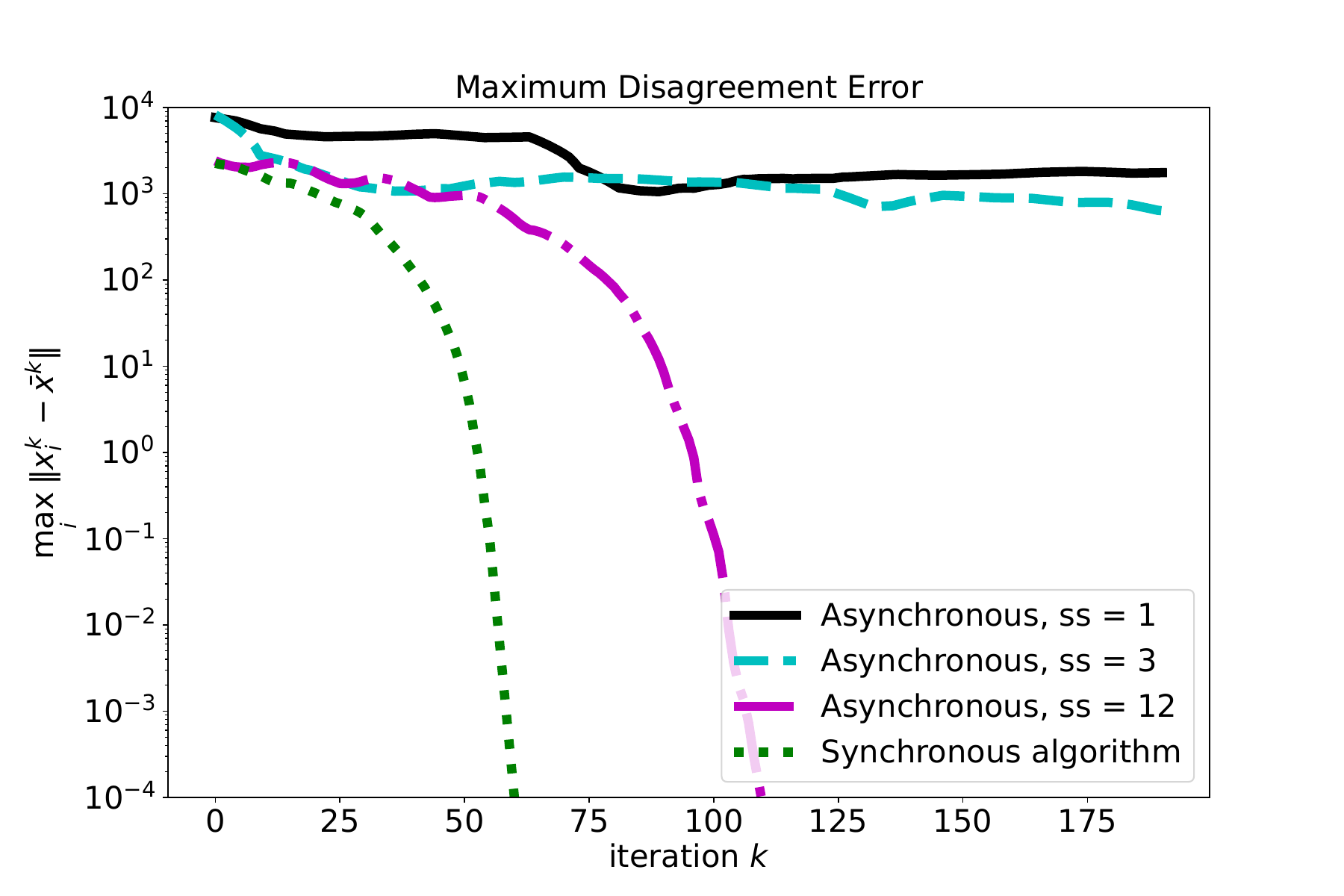}
        % \caption{$\log_{10}D^r(\vx^r)$ vs. the number of iterations}
        \label{fig:D_PCA}}
    % \end{subfigure}
  \caption{The performance of the centralised ADMM \cite{hong2017distributed} versus the proposed decentralised algorithm for solving the PCA problem for different subset sizes ($ss$). (a) $\log_{10}G^r(\bar{\vx}^r)$ vs. the number of iterations, (b) $\log_{10}D^r(\vx^r)$ vs. the number of iterations.}
  \label{fig:PCA}
\end{figure}

In our numerical experiment, we set $n=20$, $p=500$, $\lambda = 10$, and $m_i = 100$ for all $i\in [n]$. Each element of matrix $P_i$ was independently generated from a Gaussian distribution $\cN(0, 0.1^2)$. To ensure convergence, we selected the penalty parameter $\beta$ such that $\beta > 2\max_{i\in [n]} \lambda_{\text{max}}(P_i^T P_i)$. The communication graph took the form of a ring graph as in the Phase Retrieval example. As a benchmark, we also implemented the centralised algorithm from \cite{hong2017distributed}. The results are presented in Figs \ref{fig:D_PCA} and \ref{fig:G_PCA} for 50 runs of the algorithm. It can be observed that the asynchronous algorithm converges for the distributed PCA problem and the higher the number of activated agents at each iteration the faster the convergence rate and the algorithm is faster than the centralised one in \cite{hong2017distributed}. It is noteworthy that our algorithm involves more communication at each round compared to the centralised case as each node keeps a dual variable corresponding to each of its neighbours and transmits them to the neighbours.

\section{Conclusion} \label{sec: conclusion}
\noindent This paper introduces an asynchronous distributed ADMM approach for solving the non-convex distributed optimisation problem. The proposed method eliminates the need for centralised coordination, enabling any agent to perform local updates with information received from a subset of its physically connected neighbours. We establish sufficient conditions under which asynchronous ADMM asymptotically satisfies the first-order optimality conditions with probability one. Applications of the proposed approach to solve the Phase Retrieval and PCA problems demonstrate its efficiency and correctness, making it more suitable and scalable for distributed optimisation in large-scale systems. The results presented in this paper lay the theoretical foundation for studies on asynchronous distributed ADMM and its relation to the Block Coordinate DRS method.

\printbibliography

\section{Appendix} 

\subsection{Background Definitions and Lemmas} \label{sec: background}
The \emph{epigraph} of an extended real-valued function $f: \RR^n \to \overline{\RR}$ is defined as the set $\text{epi} \ f := \{(\vx, t) \in \RR^n \times \RR \ | \ f(\vx) \leq t \}$. A function $f$ is called proper if $f(\vx)<\infty$ for at least one $\vx\in \RR^n$ and $f(\vx)>-\infty$ for all $\vx\in \RR^n$. It is said to be \emph{closed} or equivalently \emph{lower semicontinuous (lsc)} if its epigraph is a closed set in $\RR^{n+1}$, see \cite[Thm. 1.6]{rockafellar1998VariationalAnalysis}. We write $f\in C^{1,1}_L$ to indicate the class of functions $f:\RR^n \to \RR$ that are differentiable and have Lipschitz continuous gradients with parameter $L$. For simplicity, we say that such an $f$ is \emph{L-smooth} or \emph{smooth}. 

We use the following important lemma on smooth functions.
\begin{lemma} \label{descent_lemma_1} 
        (Descent Lemma \cite[Proposition~A.24]{bertsekas1997nonlinear}) Let the function $f: \RR^n \to \RR$ be a $L_f$-smooth function. Then for every $\vx$, $\vy \in \RR^n$ the following holds
        \begin{align*}
            f(\vy) \leq f(\vx) + \langle \nabla f(\vx), \vy - \vx\rangle + \frac{L_f}{2} \|\vy - \vx \|^2.
        \end{align*}
\end{lemma}
The development of optimisation methods based on splitting relies significantly on the use of proximal mappings of functions. We introduce the concept and several relevant lemmas that has been employed in the analysis.
\begin{definition}
       The \emph{proximal mapping} of function $f:\RR^n \to \overline{\RR}$ with parameter $\gamma$ is a set valued mapping defined as $\prox_{\gamma f}(\vx) : = \argmin_{\vu\in \RR^n} \{f(\vu) + \frac{1}{2 \gamma} \|\vu - \vx \|^2\}$.  The function $f$ is called \emph{prox-bounded} if there exist $\gamma_f$ that $f+\frac{1}{2\gamma}\| \cdot \|^2$ is lower bounded for all $\gamma \in (0, \gamma_f)$.
\end{definition}
\begin{lemma} \label{lemma: prox properties}
    The operator $\prox_{\gamma f}$ is single valued for smooth function $f$ and it holds that $\vy = \prox_{\gamma f}(\vx)$ if and only if $\vx = \vy + \gamma \nabla f(\vy)$.
\end{lemma}
\begin{lemma} \cite[Proposition~2.3]{themelis2020douglas} \label{lemma: prox coercivity}
    Let function $f$ be a $L$-smooth function such that $f+\frac{\ell}{2}\| \cdot \|^2$ is convex. For $\vy_i \in \prox_{\gamma f} (\vx_i)$, $i=1,2$, we have
    \begin{enumerate}
        \item $\langle \vx_2 - \vx_1, \vy_2 - \vy_1 \rangle \geq (1-\gamma \ell) \|\vy_2 - \vy_1 \|^2$,
        \item $\langle \vx_2 - \vx_1, \vy_2 - \vy_1 \rangle \geq \frac{1}{1+\gamma L} \|\vx_2 - \vx_1 \|^2$.
    \end{enumerate}
\end{lemma}
\begin{lemma}\cite[Example~5.23]{rockafellar1998VariationalAnalysis} \label{lemma: osc of prox}
    For a proper, lower semicontinuous, and prox-bounded function $f:\RR^n \to \overline{\RR}$, the operator $\prox_{\gamma f}$ is outer semicontinuous, i.e. for every $\bar{\vx} \in \RR^n$ 
    \begin{align*}
        \limsup_{\vx \to \bar{\vx}} \prox_{\gamma f}(\vx)\subset \prox_{\gamma f}(\bar{\vx}).
    \end{align*}
\end{lemma}
\begin{definition}
    Given $f: \RR^n\to \overline{\RR}$ and $A\in \RR^{m\times n}$, the infimal post-composition (IPC) of $f$ by $A$, that is $A\rhd f: \RR^m \to \RR \cup \{\pm \infty\}$, is defined as
    \begin{align}
        (A \rhd f)(\vs) := \inf_{\vx \in \{\vx \mid A\vx = \vs \}} f(\vx).
    \end{align}
\end{definition}
    The following lemma on the proximal mapping of infimal post-composition is of great importance in the development of the algorithms.
    \begin{lemma} (\cite[Ch. ~3]{ryu2022large}) \label{lemma: prox of IPC}
        Let $f: \RR^n\to \overline{\RR}$, $A\in \RR^{m\times n}$, and $\vs \in \RR^m$. If the set $\argmin_{\vx \in \RR^n} f(\vx) + \frac{1}{2}\|A\vx - \vs \|^2$ is nonempty, then the following holds
        \begin{align*}
            \prox_{A \rhd f}(\vs) = A \{\argmin_{\vx \in \RR^n} {f(\vx) + \frac{1}{2}\|A\vx - \vs \|^2}\}.
        \end{align*}
    \end{lemma}

Next, we state a martingale convergence theorem that is
critical in the proof of Theorem \ref{thm: main result}. 
\begin{lemma} \label{lemma: super MG thm}(Supermartingale Convergence Theorem \cite{robbins1971convergence})
    Let $(\Omega, \cF, \PP)$ be a probability space and $\cF_1 \subset \cF_2 \subset \dots \subset \cF$ be a sequence of $\sigma$-algebras on $\cF$. Let $Y_n$, $\beta_n$ and $Z_n$ be nonnegative $\cF_n$-measurable random variables such that for each $n$ we have
    \begin{align}
        \EE(Y_{n+1} \mid \cF_n) \leq Y_n (1+\beta_n) +W_n - Z_n,
    \end{align}
    Then the followings hold with probability one on the event $\{\sum_{n\geq 1} \beta_n < \infty, \sum_{n \geq 1} W_n < \infty \}$:
    \begin{enumerate} [label=\roman*.]  
        \item the $\lim_{n\to \infty} Y_n$ exists and is finite, and 
        \item $\sum_{n\geq 1}Z_n < \infty$.
    \end{enumerate}
\end{lemma}

\end{document}